\documentclass{amsart}
\usepackage{latexsym}
\usepackage{amsmath, amssymb}
\usepackage[all]{xy}
\newtheorem{dfs}{Definition}[section]
\newtheorem{lms}{Lemma}[section]
\newtheorem{thms}{Theorem}[section]
\newtheorem{props}{Proposition}[section]

\newtheorem{cors}{Corollary}[section]

\newtheorem{uprops}{Universal Property}[section]

\begin{document}

\title{Dimension growth for $C^*$-algebras}
\author{Andrew S. Toms}
\address{Department of Mathematics and Statistics, University of New Brunswick, 
Fredericton, New Brunswick, Canada, E3B 5A3}
\email{atoms@unb.ca}

\begin{abstract}
We introduce the growth rank of a $C^*$-algebra --- a $(\mathbb{N} \cup \{\infty\})$-valued
invariant whose minimal instance
is equivalent to the condition that an algebra absorbs the Jiang-Su
algebra $\mathcal{Z}$ tensorially --- and prove that its range is exhausted by
simple, nuclear $C^*$-algebras.  As consequences we obtain a well developed
theory of dimension growth for approximately homogeneous (AH) $C^*$-algebras,
establish the existence of simple, nuclear,
and non-$\mathcal{Z}$-stable $C^*$-algebras which are not tensorially prime, 
and show the assumption of $\mathcal{Z}$-stability to be particularly natural
when seeking classification results for nuclear $C^*$-algebras
via $\mathrm{K}$-theory.

The properties of the growth rank lead us to propose a universal 
property which can be considered inside any class of unital and nuclear
$C^{*}$-algebras.  We prove that $\mathcal{Z}$ satisfies this
universal property
inside a large class of locally subhomogeneous algebras, representing
the first uniqueness theorem for $\mathcal{Z}$ 
which does not depend on the classification theory of nuclear
$C^{*}$-algebras.
\end{abstract}

\maketitle

\section{Introduction}\label{intro}

In the late 1980s, Elliott conjectured that separable nuclear
$C^*$-algebras would be classified by $\mathrm{K}$-theoretic invariants.  He 
bolstered his claim by proving that certain inductive limit $C^*$-algebras
(the $\mathrm{A}\mathbb{T}$ algebras of real rank zero, \cite{El3}) were so classified, generalising broadly
his seminal classification of approximately finite-dimensional (AF) algebras
by their scaled, ordered $\mathrm{K}_0$-groups (\cite{El1}, 1976).  
His conjecture was confirmed in the case of simple algebras
throughout the 1990s and early 2000s.  Highlights include the
Kirchberg-Phillips classification of purely infinite simple $C^*$-algebras
satisfying the Universal Coefficients Theorem (UCT), the Elliott-Gong-Li
classification of simple unital approximately homogeneous (AH) algebras of bounded topological
dimension, and Lin's classification of certain tracially AF algebras.
The classifying invariant, consisting of topological
$\mathrm{K}$-theory, traces (in the stably finite case), and
a connection between them is known as the Elliott invariant.  
(See \cite{R1} for a thorough introduction to this invariant and
the classification program for separable, nuclear $C^{*}$-algebras.)

Counterexamples to Elliott's conjecture appeared first in 2002:  R{\o}rdam's 
construction of a simple, nuclear, and separable $C^*$-algebra containing a finite
and an infinite projection (\cite{R3}) was followed by two stably finite
counterexamples, due to the author (\cite{To1}, \cite{To2}).  (The
second of these (\cite{To2}) shows that Elliott's conjecture will not 
hold even after adding to the Elliott invariant 
every continuous (with respect to direct limits) homotopy invariant
functor from the category of $C^{*}$-algebras.)  The salient common feature of these counterexamples
is their failure to absorb the Jiang-Su algebra $\mathcal{Z}$ tensorially.  
The relevance of this property to Elliott's classification program derives from the following fact:
taking the tensor product of a simple unital $C^*$-algebra $A$ with $\mathcal{Z}$
is trivial at the level of the Elliott invariant when $A$ has weakly unperforated
ordered $\mathrm{K}$-theory (\cite{GJS}), and so the Elliott conjecture predicts 
that any simple, separable, unital, and nuclear $A$ satisfying this $\mathrm{K}$-theoretic condition
will also satisfy $A \otimes \mathcal{Z} \cong A$. This last condition is known
as \emph{$\mathcal{Z}$-stability}, and any $A$ satisfying it is said to be \emph{$\mathcal{Z}$-stable}.
In recent work with Wilhelm Winter, the author has proved that every class of unital,
simple, and infinite-dimensional $C^*$-algebras for which the Elliott conjecture is
so far confirmed consists of $\mathcal{Z}$-stable algebras.  The emerging consensus, 
suggested first by R{\o}rdam and well supported by these results, is that the Elliott 
conjecture should hold for simple, nuclear, separable, and $\mathcal{Z}$-stable $C^*$-algebras.

A recurring theme in theorems confirming the Elliott conjecture is that of minimal rank.
There are various notions of rank for $C^*$-algebras --- the real rank, the stable rank,
the tracial topological rank, and the decomposition rank --- which attempt to capture 
a non-commutative version of dimension.  A natural and successful approach to proving
classification theorems for separable and nuclear $C^*$-algebras has been to assume that
one or more of these ranks is minimal (see \cite{EG} and \cite{Li1}, for instance).
But there are examples which
show these minimal rank conditions to be variously too strong or too weak to characterise
those algebras for which the Elliott conjecture will be confirmed.  One wants to assume
 $\mathcal{Z}$-stability instead, but a fair objection has been that this assumption seems
unnatural.

In the sequel, we situate $\mathcal{Z}$-stability as the minimal instance of a well-behaved
rank for $C^*$-algebras, which we term the \emph{growth rank}.  The growth rank measures
``how far'' a given algebra is from being $\mathcal{Z}$-stable, and inherits excellent
behaviour with respect to common operations from the robustness of $\mathcal{Z}$-stability.
Our terminology is motivated by the fact that the growth rank may be viewed as 
a theory of dimension growth for AH algebras, and, more generally, locally type-I
$C^*$-algebras.  We prove that for every $n \in \mathbb{N} \cup \{\infty\}$, 
there is a simple, separable, and nuclear $C^{*}$-algebra $A_{n}$ with
growth rank equal to $n$.  The algebras constructed in the proof of this theorem 
are entirely new, and rather exotic;  for all but two of them, the other ranks for
$C^*$-algebras above are simultaneously infinite.  We use these algebras to obtain 
the unexpected (see \cite{R3}):  a simple, nuclear, and non-$\mathcal{Z}$-stable 
$C^*$-algebra which is not tensorially prime. 

Motivated by the properties of the growth rank, we propose a pair of conditions on a
unital and nuclear $C^{*}$-algebra $A$ which constitute a universal
property.  The first of these conditions is known to hold for
$\mathcal{Z}$.  We verify the second condition for $\mathcal{Z}$
inside a large class of separable, unital, nuclear, and
locally subhomogeneous $C^*$-algebras which, significantly, 
contains projectionless algebras.  This represents the first 
uniqueness theorem for $\mathcal{Z}$ among projectionless algebras
which does not depend on the classification of such algebras via 
the Elliott invariant.

Our paper is organised as follows:  in section \ref{gr} we introduce the growth rank
and establish its basic properties;  in section \ref{tdgsec} we show that the growth
rank may be viewed as an abstract version of dimension growth for AH algebras;  in section 
\ref{rangesec} we establish the range of the growth rank, and consider
the growth rank of some examples; tensor factorisation and the existence of 
a simple, nuclear, and non-$\mathcal{Z}$-stable
$C^*$-algebra which is not tensorially prime are contained in section \ref{tens}; 
two universal propoerties for a simple, separable, unital, and nuclear 
$C^{*}$-algebra are discussed in section \ref{uniprops}, and the second of
these is shown to satisfied by $\mathcal{Z}$ inside certain classes of algebras;  
 connections between the
growth rank and other ranks for $C^*$-algebras are drawn in section \ref{otherranks},
and it is argued that the growth rank is connected naturally to Elliott's
conjecture.

\emph{Acknowledgements.} The author would like to thank Bruce Blackadar,
Mikael R{\o}rdam, and Wilhelm Winter for many helpful comments and suggestions.

\section{The growth rank of a $C^*$-algebra}\label{gr}

Recall that the Jiang-Su algebra $\mathcal{Z}$ is a simple, unital, nuclear, and infinite-dimensional
$C^*$-algebra which is $\mathrm{KK}$-equivalent to the
complex numbers (cf. \cite{JS1}).  We say that a $C^*$-algebra $A$ is $\mathcal{Z}$-stable
if $A \otimes \mathcal{Z} \cong A$.
The existence of simple, nuclear, separable, and non-elementary $C^*$-algebras
which are not $\mathcal{Z}$-stable was established by Villadsen 
in \cite{V1}. 

\begin{dfs}\label{growthrank}
Let $A$ be a $C^*$-algebra.  The \emph{growth rank}
$\mathrm{gr}(A)$ is the least natural number $n$ such that
\[
A^{\otimes n} \stackrel{\mathrm{def}}{=} \underbrace{A \otimes \cdots \otimes A}_{n \ \mathrm{times}}
\]
is $\mathcal{Z}$-stable, assuming the minimal tensor product.  If no such integer exists, 
then say $\mathrm{gr}(A) = \infty$.
\end{dfs}

The growth rank is most interesting for $C^*$-algebras without finite-dimensional
representations, as these are the only algebras whose finite tensor powers may be
$\mathcal{Z}$-stable.  Thus, the growth rank differs significantly from
other notions of rank for nuclear $C^*$-algebras --- the stable rank, the real rank, the tracial
topological rank, and the decomposition rank --- see \cite{BP}, \cite{Ri1}, \cite{Li3} and \cite{KW}, 
respectively, for definitions and basic properties  --- 
in that it is not proportional to the covering dimension of the spectrum in the commutative case.  
Rather, it is designed to recover information about $C^*$-algebras which are pathological with 
respect to the Elliott conjecture.     

The permanence properties of $\mathcal{Z}$-stability, most of them established in
\cite{TW1}, show the growth rank to be remarkably well behaved with respect to common operations.

\begin{thms}\label{grprop}
Let $A$, $B$ be separable, nuclear $C^*$-algebras, $I$ a closed two-sided ideal of $A$,
$H$ a hereditary subalgebra of $A$, and $k \in \mathbb{N}$.  Then,
\begin{enumerate}
\item[(i)] $\mathrm{gr}(H) \leq \mathrm{gr}(A)$;
\item[(ii)] $\mathrm{gr}(A/I) \leq \mathrm{gr}(A)$;
\item[(iii)] $\mathrm{gr}(A) = \mathrm{gr}(A \otimes \mathrm{M}_k) = \mathrm{gr}(A \otimes \mathcal {K})$;
\item[(iv)] $\mathrm{gr}(A \otimes B) \leq \mathrm{min} \{\mathrm{gr}(A),\mathrm{gr}(B)\}$;
\item[(v)] $\mathrm{gr}(A \oplus B) \leq \mathrm{gr}(A) + \mathrm{gr}(B)$;
\item[(vi)] if $A_1,\ldots,A_k$ are hereditary subalgebras of $A$, then
\[
\mathrm{gr}\left(\oplus_{i=1}^k A_i \right) \leq \mathrm{gr}(A);
\]
\item[(vii)] if $A$ is the limit of an inductive sequence $(A_i,\phi_i)$, 
where $A_i$ is separable, nuclear and satisfies $\mathrm{gr}(A_i) \leq n$ 
for each $i \in \mathbb{N}$, then $\mathrm{gr}(A) \leq n$;
\item[(viii)] if $\mathrm{gr}(I) = \mathrm{gr}(A/I) = 1$, then $\mathrm{gr}(A) = 1$.
\end{enumerate}
\end{thms}

\begin{proof} 
(i) and (ii) are clearly true if $\mathrm{gr}(A) = \infty$.  Suppose 
that $\mathrm{gr}(A) = n \in \mathbb{N}$, so that $A^{\otimes n}$ is 
$\mathcal{Z}$-stable.  Since $H^{\otimes n}$ is a hereditary subalgebra of $A^{\otimes n}$ we conclude
that it is $\mathcal{Z}$-stable by Corollary 3.3 of \cite{TW1} ---
$\mathcal{Z}$-stability passes to hereditary subalgebras.  
(ii) follows from Corollary 3.1 of \cite{TW1} after noticing that $(A/I)^{\otimes n}$
is a quotient of $A^{\otimes n}$.

(iii) is Corollary 3.2 of \cite{TW1}.

For (iv), suppose that $\mathrm{gr}(A) \leq \mathrm{gr}(B)$.
Then, 
\[
(A \otimes B)^{\otimes \mathrm{gr}(A)} \cong A^{\otimes \mathrm{gr}(A)} \otimes 
B^{\otimes \mathrm{gr}(A)}
\]
is $\mathcal{Z}$-stable since is the tensor product of two algebras, one of 
which --- $A^{\otimes \mathrm{gr}(A)}$ --- is $\mathcal{Z}$-stable.  

For (v), one can use the binomial theorem to write
\[
(A \oplus B)^{\otimes \mathrm{gr}(A) + \mathrm{gr}(B)} \cong \bigoplus_{i=0}^{\mathrm{gr}(A) + \mathrm{gr}(B)}
A^{\otimes i} \otimes B^{\otimes \mathrm{gr}(A) + \mathrm{gr}(B) - i}.
\]
For each $0 \leq i \leq \mathrm{gr}(A) + \mathrm{gr}(B)$ one has that either $i \geq \mathrm{gr}(A)$
or $\mathrm{gr}(A) + \mathrm{gr}(B)- i \geq \mathrm{gr}(B)$, whence $A^{\otimes i} 
\otimes B^{\otimes \mathrm{gr}(A) + \mathrm{gr}(B) - i}$ is $\mathcal{Z}$-stable.  It follows that
$(A \oplus B)^{\mathrm{gr}(A) + \mathrm{gr}(B)}$ is $\mathcal{Z}$-stable, as required.

For (vi) we use the fact that
\[
\left( \oplus_{i=1}^k A_i \right)^{\otimes \mathrm{gr}(A)}
\]
is a direct sum of algebras of the form
\[
A_1^{\otimes n_1} \otimes A_2^{\otimes n_2} \otimes \cdots \otimes A_k^{\otimes n_k}, \ \ \sum_{i=1}^k n_i = \mathrm{gr}(A),
\]
and each such algebra is a hereditary subalgebra of $A^{\otimes \mathrm{gr}(A)}$.  The desired 
conclusion now follows from (2).

(vii) is Corollary 3.4 of \cite{TW1}, while (viii) is Theorem 4.3 of the same paper.
\end{proof}

We defer our calculation of the range of the growth rank until section
\ref{rangesec}

\section{The growth rank as abstract dimension growth}\label{tdgsec}

In this section we couch the growth rank as a measure
of dimension growth in the setting of AH algebras.  
Recall that an unital AH algebra is an inductive limit 
\begin{equation}\label{stddecomp}
A \cong \lim_{i \to \infty}(A_i,\phi_i)
\end{equation}
where $\phi_{i}:A_{i} \to A_{i+1}$ is an unital $*$-homomorphism and
\begin{equation}\label{stdAi}
A_i := \bigoplus_{l=1}^{m_i} p_{i,l} (\mathrm{C}(X_{i,l}) \otimes
\mathcal{K}) p_{i,l}
\end{equation}
for compact connected Hausdorff spaces $X_{i,l}$ of finite covering
dimension, projections $p_{i,l} \in 
\mathrm{C}(X_{i,l}) \otimes \mathcal{K}$ ($\mathcal{K}$ is the algebra
of compact operators on a separable Hilbert space $\mathcal{H}$), and natural numbers $m_i$.
Put 
\[
\phi_{ij}:= \phi_{j-1} \circ \phi_{j-2} \circ \cdots \circ \phi_i.
\]
We refer to this
collection of objects and maps as a \emph{decomposition} for $A$.  
Decompositions for $A$ are highly non-unique.  The proof of Theorem
2.5 in \cite{Go1} shows that one may assume the $X_{i,l}$ above to be
finite CW-complexes.  We make this assumption throughout the sequel.

When we speak of dimension growth for an AH algebra we are referring,
roughly, to the asymptotic behaviour of the ratios
\[
\frac{\mathrm{dim}(X_{i,l})}{\mathrm{rank}(p_{i,l})},
\]
assuming, due to the non-unique nature of decompositions for $A$,
that we are looking at a decomposition for which these ratios grow
at a rate close to some lower limit.
If there exists a decomposition for $A$ such that
\begin{equation}\label{sdgdef}
\lim_{i \to \infty} \mathrm{max}_{1 \leq l \leq m_i}\left\{\frac{\mathrm{dim}(X_{i,l})}{\mathrm{rank(p_{i,l})}}\right\} = 0,
\end{equation}
then $A$ is said to have \emph{slow dimension growth}.  This definition appeared first in
\cite{BDR}.  As it turns out, this definition is not suitable for
non-simple algebras, at least from the point of view that slow
dimension growth should entail good behaviour in ordered
$\mathrm{K}$-theory.  This is pointed out by Goodearl in \cite{Go1}, and a second,
more technical definition of slow dimension growth is introduced.
We are interested in a demonstration of principle --- that the growth
rank yields a theory of dimension growth for AH algebras --- and so 
will limit technicalities by restricting our attention to direct 
sums of simple and unital AH algebras.
In this setting, Goodearl's definition and the one above coincide.

(Slow dimension growth or, occasionally, a slightly stronger version
thereof, is an essential hypothesis in classification theorems
for simple unital AH algebras.)

Observe that taking the tensor product of an
unital AH algebra with itself reduces dimension
growth.  Indeed, for compact Hausdorff spaces $X$ and $Y$ and natural numbers $m$ and $n$ one has
\[
\mathrm{M}_n(\mathrm{C}(X)) \otimes \mathrm{M}_m(\mathrm{C}(Y)) \cong
\mathrm{M}_{nm}(\mathrm{C}(X \times Y));
\]
the dimension of the unit in a tensor product of two homogeneous $C^*$-algebras 
is the product of the dimensions of the units, while the dimension of the spectrum
of the tensor product is the sum of the dimensions of the spectra.  If, for instance, one
has a sequence of natural numbers $n_i \stackrel{i \to \infty}{\longrightarrow} \infty$ and
an unital inductive limit algebra $A = \lim_{i \to \infty}(A_i,\phi_i)$
where $A_i = M_{n_i}(\mathrm{C}(X_i))$ and $\mathrm{dim}(X_i) = n_i^{m}$, then
$A^{\otimes m+1}$ has slow dimension growth, despite the fact that $A$ may not;
$A^{\otimes m+1}$ is an inductive limit of the 
building blocks
\[
A_i^{\otimes m+1} \cong M_{n_i^{m+1}}(\mathrm{C}((X_i)^{m+1}),
\]
and 
\[
\frac{\mathrm{dim}((X_i)^{m+1})}{n_i^{m+1}} = \frac{(m+1)n_i^{m}}{n_i^{m+1}} \stackrel{i \to \infty}{\longrightarrow} 0.
\]
We use this observation to define a concrete measure of dimension growth for unital AH algebras.

\begin{dfs}\label{tdg}
Let $A$ be an unital AH algebra.  Define the \emph{topological dimension growth} 
$\mathrm{tdg}(A)$ to be the least non-negative integer $n$ such that
$A^{\otimes n}$ has slow dimension growth, if it exists.  If no such integer exists, then say
$\mathrm{tdg}(A) = \infty$. 
\end{dfs}

Roughly, an unital AH algebra with finite topological
dimension growth $n$ has a decomposition for which
\begin{equation}\label{propor}
\mathrm{dim}(X_{i,l}) \propto \mathrm{rank}(p_{i,l})^{n-1},
\end{equation}
and no decomposition for which (\ref{propor}) holds
with $n$ replaced by $m < n$.  One might say that such an algebra
has ``polynomial dimension growth of order $n-1$''.  Similarly, an algebra
for which $\mathrm{tdg}=\infty$ has ``exponential dimension growth''.   

We now compare the properties of the topological dimension growth to 
those of the growth rank.

\begin{lms}\label{sdgsum}
Let $A$ and $B$ be unital AH algebras with slow dimension growth.  Then, $A \oplus B$ has
slow dimension growth.
\end{lms}

\begin{proof}
Straightforward.
\end{proof}

\begin{lms}\label{sdgprod}
Let $A$ and $B$ be simple and unital AH algebras, and suppose that $A$ has slow dimension
growth.  Then, $A \otimes B$ has slow dimension growth.
\end{lms}

\begin{proof}
We exploit the fact that there is considerable freedom in choosing an inductive
limit decomposition for $A \otimes B$, even after fixing decompositions for $A$ and $B$.
Let $A$ be decomposed as in (\ref{stddecomp}) and (\ref{stdAi}), and let  
\[
B \cong \lim_{j \to \infty}\left( \bigoplus_{s=1}^{n_j} q_{j,s} \mathrm{M}_{t_{j,s}}(\mathrm{C}(Y_{j,s})) q_{j,s}, \psi_j \right)
\] 
be a decomposition of $B$, where the $Y_{j,s}$ are connected compact Hausdorff spaces and
the $q_{j,s} \in \mathrm{M}_{t_{j,s}}(\mathrm{C}(Y_{j,s}))$ are projections.
Put
\[
B_j:=\bigoplus_{s=1}^{n_j} q_{j,s} \mathrm{M}_{t_{j,s}}(\mathrm{C}(Y_{j,s})) q_{j,s}.
\]
For any strictly increasing sequence $(r_i)$ of natural numbers one has
\[
A \otimes B \cong \lim_{i \to \infty}(A_{r_i} \otimes B_i,\phi_{r_i,r_{i+1}} \otimes \psi_i).
\]
Put 
\[
M_i \stackrel{\mathrm{def}}{=} \mathrm{max}_{1 \leq s \leq n_i} \{\mathrm{dim}(Y_{i,s})\}.
\]

The simplicity of $A$ implies that for any $N \in \mathbb{N}$, there exists $i_N \in \mathbb{N}$
such that $\mathrm{rank}(p_{i,l}) \geq N$, $\forall i \geq i_N$.  Choose the sequence $(r_i)$ so
that  
\[
\mathrm{min}_{1 \leq l \leq m_{r_i}} \{\mathrm{dim}(p_{r_i,l})\} \geq 2^i M_i.
\]
A typical direct summand of $A_{r_i} \otimes B_i$ with connected spectrum has the form
\[
(p_{r_i,l} \otimes q_{i,s})\left(\mathrm{M}_{k_{r_i,l}t_{i,s}}(\mathrm{C}(X_{r_i,l} \times Y_{i,s}))\right)
(p_{r_i,l} \otimes q_{i,s}),
\]
whence the condition that $A \otimes B$ has slow dimension growth amounts to the condition that
\[
\liminf_{i \to \infty} \mathrm{max}_{l,s} \left\{ \frac{\mathrm{dim}(X_{r_i,l}) + \mathrm{dim}(Y_{i,s})}
{\mathrm{rank}(p_{r_i,l}) \mathrm{rank}(q_{i,s})} \right\} = 0.
\]
We have that
\[
\mathrm{max}_{l,s} \left\{ \frac{\mathrm{dim}(X_{r_i,l}) + \mathrm{dim}(Y_{i,s})}
{\mathrm{rank}(p_{r_i,l}) \mathrm{rank}(q_{i,s})} \right\}
\]
is dominated by 
\[
\mathrm{max}_{l} \left\{ \frac{\mathrm{dim}(X_{r_i,l})}
{\mathrm{rank}(p_{r_i,l}) \mathrm{rank}(q_{i,s})} \right\}+
\mathrm{max}_{s} \left\{ \frac{\mathrm{dim}(Y_{i,s})}
{\mathrm{rank}(p_{r_i,l}) \mathrm{rank}(q_{i,s})} \right\}.
\]
In the above sum the first term tends to zero by virtue of $A$ having slow dimension growth, 
while the second tends to zero by our choice of $(r_i)$.  We conclude that $A \otimes B$
has slow dimension growth, as desired.
\end{proof}

\begin{thms}\label{tdgprop}
Let $A$, $B$ be simple and unital AH algebras.  Then,
\begin{enumerate}
\item[(i)] $\mathrm{tdg}(A \otimes B) \leq \mathrm{min}
\{\mathrm{tdg}(A), \mathrm{tdg}(B)\}$;
\item[(ii)] $\mathrm{tdg}(A \oplus B) \leq \mathrm{tdg}(A) +
\mathrm{tdg}(B)$.
\end{enumerate}
\end{thms}

\begin{proof}

For (i), suppose that $\mathrm{tdg}(A) \leq \mathrm{tdg}(B)$.  Then,
\[
(A \otimes B)^{\otimes \mathrm{tdg}(A)} \cong (A^{\otimes \mathrm{tdg}(A)}) \otimes (B^{\otimes \mathrm{tdg}(A)}).
\]
Since $A^{\otimes \mathrm{tdg}(A)}$ has slow dimension growth by defintion, the 
right hand side of the equation above has slow dimension growth by Lemma \ref{sdgprod}.

For (ii), use the binomial theorem to write
\[
(A \oplus B)^{\otimes \mathrm{tdg}(A) + \mathrm{tdg}(B)} \cong \bigoplus_{i=0}^{\mathrm{tdg}(A) + \mathrm{tdg}(B)}
A^{\otimes i} \otimes B^{\otimes \mathrm{tdg}(A) + \mathrm{tdg}(B)- i}.
\]
Notice that each direct summand of the right hand side above has slow dimension growth
by part $(2)$ of this proposition, whence the entire direct sum has slow dimension growth by
Lemma \ref{sdgsum}.
\end{proof}

As far as direct sums of simple unital AH algebras are concerned, the properties of the growth rank
agree with those of the topological dimension growth, despite that fact that $\mathcal{Z}$-stability
and slow dimension growth are not yet known to be equivalent for simple, unital and
infinite-dimensional AH algebras. 

Next, we prove that the topological dimension growth and the growth rank 
often agree.  Recall that a Bauer simple is a compact metrizable Choquet simplex $S$ whose 
extreme boundary $\partial_eS$ is compact.  The set $Aff(S)$ of continuous affine real-valued functions on
$S$ are in bijective correspondence with continuous real-valued functions on $\partial_eS$.  
The bijection is given by the map which assigns to a 
continuous affine function $f$ on $S$, the continuous function $\hat{f}:\partial_eS \to \mathbb{R}$
given by 
\[
\hat{f}(\tau) = f(\tau), \ \forall \tau \in \partial_eS.
\]

\begin{props}\label{dg-greaterthan-gr}
Let $A$ be a simple, unital and infinite-dimensional AH algebra.  Suppose that the
simplex of tracial states $\mathrm{T}A$ is a Bauer simplex, and that
the image of $\mathrm{K}_0(A)$ in $\mathrm{C}_{\mathbb{R}}(\partial_e \mathrm{T}A)$
is uniformly dense.  Then, 
\[
\mathrm{tdg}(A) = \mathrm{gr}(A).
\]
\end{props}

\begin{proof}
It is well known that 
\[
\partial_e \mathrm{T} A^{\otimes n} \cong \times_{i=1}^n \partial_e \mathrm{T}A,
\]
whence,
\[
\mathrm{C}_{\mathbb{R}}(\partial_e \mathrm{T} A^{\otimes n}) \cong 
\mathrm{C}_{\mathbb{R}}(\partial_e \mathrm{T}A)^{\otimes n}.
\]
Let $f_1,\ldots,f_n \in \mathrm{C}_{\mathbb{R}}(\partial_e \mathrm{T}A)$ 
be the images of elements $x_1,\ldots,x_n \in \mathrm{K}_0(A)$, respectively.
Write $x_i = [p_i]-[q_i]$ for projections $p_i,q_i \in \mathrm{M}_{\infty}(A)$,
$1 \leq i \leq n$.  Let $g_i,h_i \in \mathrm{C}_{\mathbb{R}}(\partial_e \mathrm{T}A)$
be the images of $p_i,q_i$, respectively.  Now
\begin{equation}\label{afftprod}
f_1 \otimes \cdots \otimes f_n  =  \bigotimes_{i=1}^n (g_i - h_i),  
\end{equation}
and the right hand side of the equation is a sum of elementary
tensor $\pm(r_1 \otimes \cdots \otimes r_n)$, where $r_i \in \{h_i,g_i\}$.
There are thus projections $t_i \in \{p_i,q_i\}$ such that the image of 
\[
[t_1 \otimes \cdots \otimes t_n] \in \mathrm{K}_0(A^{\otimes n})
\]
is 
\[
r_1 \otimes \cdots \otimes r_n \in \mathrm{C}_{\mathbb{R}}(\partial_e \mathrm{T}A)^{\otimes n}.
\]
Thus, the right hand side of (\ref{afftprod}) can be obtained as the image
of some $x \in \mathrm{K}_0(A^{\otimes n})$.  Given $\epsilon > 0$ and an
elementary tensor
\[
m_1 \otimes \cdots \otimes m_n \in \mathrm{C}_{\mathbb{R}}(\partial_e \mathrm{T}A)^{\otimes n} \cong
\mathrm{C}_{\mathbb{R}}(\partial_e \mathrm{T} A^{\otimes n}), 
\]
we may, by the density of the image of $\mathrm{K}_0(A)$ in 
$\mathrm{C}_{\mathbb{R}}(\partial_e \mathrm{T}A)$, choose $f_1,\ldots,f_n
\in \mathrm{C}_{\mathbb{R}}(\partial_e \mathrm{T}A)$ to satisfy
\[
|(m_1 \otimes \cdots \otimes m_n) - (f_1 \otimes \cdots \otimes f_n)| < \epsilon.
\]
It follows that the image of $\mathrm{K}_0(A^{\otimes n})$ is dense in
$\mathrm{C}_{\mathbb{R}}(\partial_e \mathrm{T} A^{\otimes n})$.
Theorem 3.13 of \cite{TW2} shows that $\mathcal{Z}$-stability and slow dimension
growth are equivalent for $A^{\otimes n}$, and the proposition follows.
\end{proof}

Following \cite{TW2}, we may drop the condition that the image of $\mathrm{K}_0$ in
$Aff(\mathrm{T}(A))$ be dense whenever $A$ has a unique tracial state.  Note that 
an algebra satisfying the hypotheses of Proposition \ref{dg-greaterthan-gr} 
need not have real rank zero, even in the case of a unique tracial state (cf. \cite{V2}). 

As the growth rank and the topological dimension growth often (probably always, in the simple
and unital case) agree, we suggest simply using the growth rank as a theory of unbounded dimension
growth for AH algebras.  It has the advantage of avoiding highly technical definitions
involving arbitrary decompositions for a given AH algebra, and works equally well for
non-simple and non-unital algebras.  

There are definitions of slow dimension growth for more general locally type-I $C^*$-algebras
--- direct limits of recursive subhomogeneous algebras (\cite{P3}), for instance --- but these
are even more technical than the definition for non-simple AH algebras.  The growth rank
seems the logical choice for defining dimension growth in these situations, too.

\section{A range result}\label{rangesec}

\begin{thms}\label{range}
Let $n \in \mathbb{N} \cup \{\infty\}$.  There exists a simple, nuclear, and non-type-I
$C^*$-algebra $A$ such that
\[
\mathrm{tdg}(A) = \mathrm{gr}(A) = n.
\]
\end{thms}

To prepare the proof of Theorem \ref{range}, we collect some basic facts about
the Euler and Chern classes of a complex vector bundle, and recall results of R{\o}rdam
and Villadsen.

Let $X$ be a connected topological space, and let $\omega$ and $\gamma$ be 
(complex) vector bundles over $X$ of fibre dimensions $k$ and $m$, respectively.
Recall that the Euler class $e(\omega)$ is an element of $H^{2k}(X;\mathbb{Z})$ with the following
properties:
\begin{enumerate}
\item[(i)] $e(\omega \oplus \gamma) = e(\omega) \cdot e(\gamma)$, where ``$\cdot$'' 
denotes the cup product in $H^{*}(X;\mathbb{Z})$;
\item[(ii)] $e(\theta_l) = 0$, where $\theta_l$ denotes the trivial vector bundle over $X$ 
of (complex) fibre dimension $l$.
\end{enumerate}
The Chern class $c(\omega) \in H^*(X;\mathbb{Z})$ is a sum
\[
c(\omega) = 1 + c_1(\omega) + c_2(\omega) + \cdots + c_k(\omega),
\]
where $c_i(\omega) \in H^{2i}(X;\mathbb{Z})$.  Its properties are similar to those of the
Euler class:
\begin{enumerate}
\item[(i)] $c(\omega \oplus \gamma) = c(\omega) \cdot c(\gamma)$;
\item[(ii)] $c(\theta_l) =1$.
\end{enumerate}
The key connection between these characteristic classes is this:  $e(\eta) =
c_1(\eta)$ for every line bundle. 

The next lemma is due essentially to Villadsen (cf. Lemma 1 of \cite{V1}), but our
version is more general.

\begin{lms}\label{vil} 
Let $X$ be a finite CW-complex, and let $\eta_1,\eta_2,\ldots,\eta_k$ be
complex line bundles over $X$.  If $l < k$ and $\prod_{i=1}^k e(\eta_i) \neq 0$, then
\[
[\eta_1 \oplus \eta_2 \oplus \cdots \oplus \eta_k] - [\theta_l] \notin \mathrm{K}^0(X)^+.
\]
\end{lms}

\begin{proof}
If $[\eta_1 \oplus \eta_2 \oplus \cdots \oplus \eta_k] - [\theta_l] \in \mathrm{K}^0(X)^+$,
then there is a vector bundle $\gamma$ of dimension $k-l$ and $d \in \mathbb{N}$ such
that
\[
\eta_1 \oplus \eta_2 \oplus \cdots \oplus \eta_k \oplus \theta_d \cong \gamma \oplus \theta_{d + k - l}.
\]
Applying the Chern class to both sides of this equation we obtain
\[
\prod_{i=1}^k c(\eta_i) = c(\gamma).
\]
Expanding the left hand side yields
\[
\prod_{i=1}^k (1 + c_1(\eta_i)) = \prod_{i=1}^k (1 + e(\eta_i)).
\]
The last product has only one term in $H^{2k}(X;\mathbb{Z})$, namely,
$\prod_{i=1}^k e(\eta_i)$, and this, in turn, is non-zero.  On the other
hand, $c(\gamma)$ has no non-zero term in $H^{2i}(X;\mathbb{Z})$ for $i > k-l$.
Thus, we have a contradiction, and must conclude that 
\[
[\eta_1 \oplus \eta_2 \oplus \cdots \oplus \eta_k] - [\theta_l] \notin \mathrm{K}^0(X)^+.
\]
\end{proof}

Let $\xi$ be any line bundle over $\mathrm{S}^2$ with non-zero Euler class --- the Hopf 
line bundle, for instance.
We recall some notation and a proposition from \cite{R3}.  For each natural number
$s$ and for each non-empty finite set 
\[
I= \{s_1,\ldots,s_k\} \subseteq \mathbb{N}
\]
define bundles $\xi_s$ and $\xi_I$ over ${\mathrm{S}^2}^m$ (for all $m \geq s$ or $m \geq \mathrm{max}
\{s_1,\ldots,s_k\}$, as appropriate) by 
\[
\xi_s = \pi_s^*(\xi), \ \ \ \ \ \xi_I = \xi_{s_1} \otimes \cdots \otimes \xi_{s_k},
\]
where $\pi_s:{\mathrm{S}^2}^{\times m} \to \mathrm{S}^2$ is the $s^{\mathrm{th}}$
co-ordinate projection. 

\begin{props}[R{\o}rdam, Proposition 3.2, \cite{R3}]\label{marriage}
Let $I_1,\ldots,I_m \subseteq \mathbb{N}$ be finite sets.  The following are equivalent:
\begin{enumerate}
\item[(i)] $e(\xi_{I_1} \oplus \xi_{I_2} \oplus \cdots \oplus \xi_{I_m}) \neq 0$.
\item[(ii)] For all subsets $F$ of $\{1,2,\ldots,m\}$ we have $| \cup_{j \in F} I_j | \geq |F|$.
\end{enumerate}
\end{props}

\noindent
(In fact, there is a third equivalence in Proposition 3.2 of \cite{R3}.  We do not require
it, and so omit it.)

\begin{proof}(Theorem \ref{range})  
The case where $n=1$ is straightforward:  any UHF algebra $\mathcal{U}$ is $\mathcal{Z}$-stable 
by the classification theorem of \cite{JS1} (or, alternatively, by Theorem 2.3 of \cite{TW1}), 
and has $\mathrm{tdg}(A)=1$ (as does any AF algebra).  

Let $ 1 \leq  n \in \mathbb{N} \cup \{\infty\}$ be given.   
We will construct an simple, unital, and infinite-dimensional AH algebra
\[
A = \lim_{i \to \infty}(A_i,\phi_i)
\]
 along the 
lines of the construction of \cite{V2}, and prove that $\mathrm{gr}(A) = \mathrm{tdg}(A) = n+1$.  
Our strategy is to prove that $A^{\otimes n}$
has a perforated ordered $\mathrm{K}_0$-group and is hence not $\mathcal{Z}$-stable by Theorem 1 of \cite{GJS}, 
while $A^{\otimes n+1}$ is tracially AF and hence approximately divisible by \cite{EGL2} and
$\mathcal{Z}$-stable by Theorem 2.3 of \cite{TW2}.  $A$ will be constructed so as to have a unique trace,
whence $\mathrm{tdg}(A) = \mathrm{gr}(A)$ by Proposition \ref{dg-greaterthan-gr}.

Let $X_1 = (\mathrm{S}^2)^{n_1}$, and, for each $i \in \mathbb{N}$, let $X_i = 
(X_{i-1})^{n_i}$, where the $n_i$ are natural numbers to be specified.  Set
\[
N_i:=\prod_{j=1}^{i} n_j
\] 
and 
\[
I_{l}^{i} := \{(l-1)N_{i}+N_{i-1}+1,\ldots,l N_i\} \subseteq \{1,\ldots,nN_i\}, \ \ l \in \{1,\ldots,n\}.
\]
We will take
\[
A_i = p_i (\mathrm{C}(X_i) \otimes \mathcal{K}) p_i
\]
for some projection $p_i \in \mathrm{C}(X_i) \otimes \mathcal{K}$ to be specified. 
The maps 
\[
\phi_i:A_i \to A_{i+1}
\] 
are constructed inductively as follows:  suppose that
$p_1,\ldots,p_i$ have been chosen, and define a map
\[
\tilde{\phi_i}:A_i \to \mathrm{C}(X_{i+1}) \otimes \mathcal{K}
\]
by taking the direct sum of the map $\gamma_i:A_i \to \mathrm{C}(X_{i+1}) \otimes \mathcal{K}$
given by
\[
\gamma_i(f)(x) = f(\omega_i(x)) 
\]
($\omega_i:X_{i+1} \to X_i$ is projection onto the first factor of $X_{i+1} = (X_i)^{n_{i+1}}$) and
$m_{i+1}$ copies of the map $\eta_i:A_i \to \mathrm{C}(X_{i+1}) \otimes \mathcal{K}$
given by
\[
\eta_i(f)(x) = f(x_i) \cdot \xi_{I_1^{i+1}}
\]
($x_i \in X_i$ is a point to be specified, and $m_{i+1}$ is
a natural number to be specified);  set 
\[
A_{i+1} := \tilde{\phi_i}(p_i) \mathrm{C}(X_{i+1}) \otimes \mathcal{K} \tilde{\phi_i}(p_i),
\]  
and let $\phi_i$ be the restriction of $\tilde{\phi_i}$ to $A_{i+1}$.  In \cite{V2} it 
is shown that by replacing the $x_i$ with various other points from $X_{i+1}$ in
a suitable manner, one can ensure a simple limit $A := \lim_{i \to \infty} (A_i,\phi_i)$.
$A$ is unital by construction.

Let $p_1$ be the projection over $X_1$ corresponding to the Whitney sum
\[
\theta_1 \oplus \xi_{I_1^1} \oplus \xi_{I_1^1}.
\]
By Proposition 3.2 of \cite{R3} (Proposition \ref{marriage} of this section) we
have that the Euler class of $\oplus_{j=1}^k \xi_{I_1^1}$ is non-zero 
whenever $k \leq n_1$.  By Lemma \ref{vil}, we have 
\[
2[\xi_{I_1^1}]-[\theta_1] \notin {\mathrm{K}_0(A_1)}^+.
\]

For each $i \in \mathbb{N}$ one has
\[
\mathrm{K}_0(A_i^{\otimes n}) \cong \mathrm{K}^0(X_i^n) \cong \mathrm{K}^0(X_i)^{\otimes n}
\]
--- the last isomorphism follows from the K{\"u}nneth formula and the
fact that $\mathrm{K}^0(X_i)$ is torsion free.  We will prove that
\begin{equation}\label{perfel}
\left( 2[\phi_{1i}(\xi_{I_1^1})] - [\phi_{1i}(\theta_1)] \right) \otimes [p_i] \otimes
[p_i] \otimes \cdots \otimes [p_i] \notin \mathrm{K}_0(A_i^{\otimes n})^+, \ \forall i \in \mathbb{N},
\end{equation}
whence $\mathrm{K}_0(A^{\otimes n})$ is a perforated ordered group, and $A^{\otimes n}$
is not $\mathcal{Z}$-stable.

Let $Y$ and $Z$ be topological spaces, and let $\eta$ and
$\beta$ be vector bundles over $Y$ and $Z$, respectively.  Let
\[
\pi_Y: Y \times Z \to Y; \ \ \ \pi_Z: Y \times Z \to Z
\]
be the co-ordinate projections, and $\pi_Y^*(\eta)$ and $\pi_Z^*(\beta)$ the
induced bundles over $Y \times Z$.  The external tensor product $\eta \hat{\otimes}
\beta$ is defined to be the internal (fibre-wise) tensor product $\pi_Y^*(\eta) \otimes \pi_Z^*(\beta)$. 
 Let $\pi_l:(X_i)^l \to X_i$ be the
$l^{\mathrm{th}}$ co-ordinate projection, and set $p_i^l := \pi_l^*(p_i)$.
The tensor product of group elements in (\ref{perfel}) corresponds to the external tensor
product of the corresponding (formal difference of) vector bundles.   In other words, proving that 
that (\ref{perfel}) holds thus amounts to proving that
\[
2[ \phi_{1i}(\xi_{I_1^1}) \otimes p_i^2 \otimes \cdots \otimes p_i^n] - [\phi_{1i}(\theta_1) \otimes p_i^2 
\otimes \cdots \otimes p_i^n] \notin \mathrm{K}^0(X_i^n)^+.
\]
Straightforward calculation shows that $\theta_1$ is a direct summand of
$\phi_{1i}(\theta_1) \otimes p_i^2 \otimes \cdots \otimes p_i^n$, for all $i \in \mathbb{N}$.
Thus,
\begin{eqnarray*}
 & & 2[ \phi_{1i}(\xi_{I_1^1}) \otimes p_i^2 \otimes \cdots \otimes p_i^n] - [\phi_{1i}(\theta_1) \otimes p_i^2 
\otimes \cdots \otimes p_i^n] \\
& \leq & 2[ \phi_{1i}(\xi_{I_1^1}) \otimes p_i^2 \otimes \cdots \otimes p_i^n] - [\theta_1],
\end{eqnarray*}
and (\ref{perfel}) will hold if
\begin{equation}\label{suffice}
 2[ \phi_{1i}(\xi_{I_1^1}) \otimes p_i^2 \otimes \cdots \otimes p_i^n] - [\theta_1] \notin \mathrm{K}^0(X_i^n)^+.
\end{equation}

We prove that (\ref{suffice}) holds by induction.  Assume that $i = 1$.
The projection $p_1^l$ corresponds
to the vector bundle $\xi_{I_l^1} \oplus \xi_{I_l^1} \oplus \theta_1 = 2 \xi_{I_l^1} \oplus \theta_1$ 
over $X_i^n \cong (\mathrm{S}^2)^{ n_1n}$. 
Now
\[
\begin{array}{rl}
\left[ \bigotimes_{l=2}^{n} p_1^l \right] & = \left[\bigotimes_{l=2}^{n} (2\xi_{I_l^1} \oplus \theta_1)\right]\\
& = \left[\bigoplus_{\emptyset \neq J \subseteq \{2,\ldots,n\}} \left(\bigotimes_{l \in J} 2\xi_{I_l^1}\right) \oplus \theta_1\right],
\end{array}
\]  
so that
\[
 2[ \phi_{1i}(\xi_{I_1^1}) \otimes p_i^2 \otimes \cdots \otimes p_i^n] = 
\left[ 2\xi_{I_1^1} \otimes \left(\bigoplus_{\emptyset \neq J \subseteq \{2,\ldots,n\}} 
\left(\bigotimes_{l \in J} 2\xi_{I_l^1}\right) \right)\right] +[2\xi_{I_1^1}].
\]
By Lemma \ref{vil} it will suffice to show that
\[
e\left(\bigoplus_{\emptyset \neq J \subseteq \{2,\ldots,n\}} 
2\xi_{I_1^1} \otimes \left(\bigotimes_{l \in J}  2\xi_{I_l^1}\right) \oplus 2\xi_{I_1^1} \right)  \neq 0.
\]
Letting $I_J$ denote the union $\cup_{l \in J} I_l^1$ we have
\begin{equation}\label{stepone}
\bigoplus_{\emptyset \neq J \subseteq \{2,\ldots,n\}} 
2\xi_{I_1^1} \otimes \left(\bigotimes_{l \in J}  2\xi_{I_l^1}\right) \oplus 2\xi_{I_1^1}
= 2\xi_{I_1^1} \oplus \bigoplus_{\emptyset \neq J \subseteq \{2,\ldots,n\}} 2^{|J|+1} \xi_{I_1^1 \cup I_J}. 
\end{equation}
We wish to apply Proposition 3.2 of \cite{R3} to conclude that the Euler class of the bundle
above is non-zero.  This will, of course, require that $n_1$ be sufficiently large.  One easily sees that 
the dimension of the bundle above is $2\cdot 3^{n-1}$.
Let $R_1 = I_1^1$,
and define a list of subsets $R_2,\ldots,R_{2 \cdot 3^{n-1}}$ of $\mathbb{N}$ by including, for each $J \subseteq 
\{2,\ldots,n\}$, $2^{|J|+1}$ copies of $I_1^1 \cup I_J$ among the $R_j$.  The $R_j$s  are the index sets of the
tensor products of Hopf line bundles appearing as direct summands in  (\ref{stepone}).
We must choose $n_1$ large enough so that, for any finite subset $F$ of $\{1,\ldots,2 \cdot 3^{n-1}\}$, we
have $|\cup_{j \in F} R_j| \geq |F|$.  Clearly, setting $n_1 = 3^n$ will suffice.  This establishes the
base case of our induction argument.

We proceed to the induction step.  By Lemma 4.1 it will suffice
to prove that 
\[
e\left(2 \phi_{1i}(\xi_{I_1^1}) \otimes \bigotimes_{l=2}^{n} p_i^l \right) \neq 0.
\]
Suppose that for all $k < i$, $n_k$ has been chosen large enough that
\[
e\left(2 \phi_{1k}(\xi_{I_1^1}) \otimes \bigotimes_{l=2}^{n} p_k^l \right) \neq 0.
\]
Put $\omega_{i,l} := \omega_i \circ \pi_l$.
By construction we have
\[
p_i^l = \omega_{i,l}^*(p_{i-1}^l) \oplus m_i \cdot \mathrm{dim}(p_{i-1}^l) \otimes \xi_{I_l^i}
\]
and
\[
\phi_{1i}(\xi_{I_1^1}) = \omega_{i,1}^*\left(\phi_{1, i-1}(\xi_{I_1^1})\right) \oplus 
m_i \cdot \mathrm{dim}\left(\phi_{1, i-1}(\xi_{I_1^1}) \right) \otimes \xi_{I_1^i}.
\]
It follows that
\begin{eqnarray*}
2\phi_{1i}(\xi_{I_1^1}) \otimes \bigotimes_{l=2}^{n} p_i^l &
= & 2\left(\omega_i^*\left(\phi_{1, i-1}(\xi_{I_1^1})\right) \oplus 
m_i \cdot \mathrm{dim}\left(\phi_{1, i-1}(\xi_{I_1^1}) \right) \otimes \xi_{I_1^i}\right) \\
& & \otimes\left( \bigotimes_{l=2}^n \left( \omega_i^*(p_{i-1}^l) \oplus m_i \cdot \mathrm{dim}(p_{i-1}^l) \otimes \xi_{I_l^i} \right) \right) \\
& = &  \left( \omega_{i,1}^*(2\phi_{1,i-1}(\xi_{I_1^1})) \otimes \bigotimes_{l=2}^{n} \omega_{i,l}^*(p_{i-1}^l) \right) \oplus B \\
& = &  \Gamma_i^*  \left( 2\phi_{1,i-1}(\xi_{I_1^1}) \otimes \bigotimes_{l=2}^{n} p_{i-1}^l \right) \oplus B,
\end{eqnarray*}
where
\[
\Gamma_i := \omega_{i,1} \times \omega_{i,2} \times \cdots \times \omega_{i,n}: (X_i)^n \to (X_{i-1})^n
\]
and $B$ is a sum of line bundles, each of which has $\xi_{I_l^i}$ as a tensor factor for some
$l \in \{1,\ldots,n\}$.  The index sets of the line bundles making up
\[
\Gamma_i^*  \left( 2\phi_{1,i-1}(\xi_{I_1^1}) \otimes \bigotimes_{l=2}^{n} p_{i-1}^l \right)
\]
are disjoint from each $I_l^i$ by construction, so by Proposition 3.2 of \cite{R3} we have
\[
e\left( \Gamma_i^*  \left( 2\phi_{1,i-1}(\xi_{I_1^1}) \otimes \bigotimes_{l=2}^{n} p_{i-1}^l \right) \oplus B\right) \neq 0
\] 
if
\[
e\left( \Gamma_i^*  \left( 2\phi_{1,i-1}(\xi_{I_1^1}) \otimes \bigotimes_{l=2}^{n} p_{i-1}^l \right) \right) \neq 0; \ \ \ 
e(B) \neq 0.
\] 
The first inequality follows from our induction hypothesis.  For the second inequality, we have
\[
\mathrm{dim}(B) < \mathrm{dim} \left( 2 \phi_{1i}(\xi_{I_1^1}) \otimes \bigotimes_{l=2}^{n} p_i^l \right) <
\mathrm{dim}(p_i^l)^n.
\]
Choosing $n_i$ \emph{just large enough} (for reasons to be made clear shortly) to ensure that
$|I_l^i| \geq {\mathrm{dim}(p_i^l)}^n$,
we may conclude by Proposition 3.2 of \cite{R3} that $e(B) \neq 0$, as desired.

The fact that $A^{\otimes n}$ has a perforated ordered $\mathrm{K}_0$-group implies
that $\mathrm{gr}(A) > n$.  We now show that $\mathrm{gr}(A) \leq n+1$.
First, we compute an upper bound on the dimension of $X_i$. 
We have chosen $n_i$ to be just large enough to ensure that $|I_1^i| \geq \mathrm{dim}(p_i)^n$.
Using the fact that $|I_1^i| = N_i - N_{i-1}$ we have  
\[
N_i \leq \mathrm{dim}(p_i)^n + 2N_{i-1}, \ \ i \in \mathbb{N},
\]
since one could otherwise reduce the size of $n_i$ by one or more.
Set $d_i := \mathrm{dim}(p_i)$ for brevity.  
$A^{\otimes n+1}$ will have slow dimension growth if
\[
\frac{(n+1)N_i}{d_i^{n+1}} \stackrel{i \to \infty}{\to} 0.
\]
We have
\[
\frac{(n+1)N_i}{d_i^{n+1}} \leq \frac{(n+1)(d_i^n + 2N_{i-1})}{d_i^{n+1}} = \frac{n+1}{d_i} + \frac{2N_{i-1}}{d_i^{n+1}},
\]
so that we need only show
\begin{equation}\label{tozero}
\frac{2N_{i-1}}{d_i^{n+1}} \stackrel{i \to \infty}{\to} 0.
\end{equation}
But $N_{i-1}$ does not depend on $m_i$, so we may make $d_i$ large enough for (\ref{tozero}) to
hold (remember that $m_i$ may be chosen before $n_i$).  Thus, $A^{\otimes n+1}$ has slow dimension growth.
$A$ has a unique tracial state by the arguments of \cite{V2}, whence so does $A^{\otimes n+1}$.  It
follows that $A^{\otimes n+1}$ is of real rank zero by the main theorem of 
\cite{BBEK}.  The reduction theorem of \cite{D} together with the classification theorem of \cite{EG} then show 
that $A^{\otimes n+1}$ is approximately divisible, whence $A^{\otimes n+1}$ is $\mathcal{Z}$-stable by Theorem 2.3 of \cite{TW2}
and $\mathrm{gr}(A) = n+1$.  Since $A$ has a unique trace, it satisfies the hypotheses 
of Proposition \ref{dg-greaterthan-gr}, whence $\mathrm{tdg}(A) = \mathrm{gr}(A) = n+1$.  
This proves Theorem \ref{range} for $n$ finite.

To produce an algebra with infinite growth rank, we follow the construction
above, but choose the $n_i$ larger at each stage.  Begin as above with the same choice of $A_1$.
Notice that the arguments above not only show that one can choose $n_i$ large enough so that
\[
 2[ \phi_{1i}(\xi_{I_1^1}) \otimes p_i^2 \otimes \cdots \otimes p_i^n] - [\theta_1] \notin \mathrm{K}^0(X_i^n)^+,
\]
but also large enough so that
\[
 2[ \phi_{1i}(\xi_{I_1^1}) \otimes p_i^2 \otimes \cdots \otimes p_i^{in}] - [\theta_1] \notin \mathrm{K}^0(X_i^{in})^+.
\]
With the latter choice of $n_i$, one has that
$\mathrm{K}_0(A^{\otimes in})$ is a perforated ordered group for every natural number $i$.  It follows that
$\mathrm{gr}(A) = \infty$.  Now Proposition \ref{dg-greaterthan-gr} shows that $\mathrm{tdg}(A) = \infty$, proving
the theorem in full. 

Finally, in the case where $\mathrm{gr}(A) \geq 3$, we modify the base spaces $X_i$ 
to facilitate stable and real rank calculations in the sequel.  Replace
$X_i$ with $X_i^{'} := X_i \times \mathbb{D}^{id_i^2}$, where $\mathbb{D}$ denotes the closed unit 
disc in $\mathbb{C}$, and replace the eigenvalue map $\omega_i$ with the map
with a map $\omega_i^{'}:X_{i+1}^{'} \to X_i^{'}$ given by the Cartesian product 
of $\omega_i$ and any co-ordinate projection $\lambda_i:\mathbb{D}^{(i+1)d_{i+1}^2} 
\to \mathbb{D}^{id_i^2}$.  On the one hand, these modifications are trivial at
the level of $\mathrm{K}_0$, whence the proof of the lower bound on the growth 
rank of $A$ carries over to our new algebra.  On the other hand, our new algebra has
\[
\frac{\mathrm{dim}(X_i^{'})}{d_i^{n+1}} = \frac{\mathrm{dim}(X_i) + id_i^2}{d_i^{n+1}}
\stackrel{i \to \infty}{\to} 0,
\]
since $n \geq 2$.  Thus, the specified adjustment to the construction of $A$ does not
increase the topological dimension growth.   
\end{proof}

We now consider the growth rank of some examples.

\hspace{3mm}

\noindent
{\bf gr = 1.}
Let $A$ be a separable, unital, and approximately divisible $C^*$-algebra.  Then
$\mathrm{gr}(A) = 1$ by Theorem 2.3 of \cite{TW2}.

Let $A$ be a simple nuclear $C^*$-algebra which is neither finite-dimensional
nor isomorphic to the compact operators.  Suppose further that $A$ belongs to
a class of $C^*$-algebras for which Elliott's classification conjecture is
currently confirmed (cf. \cite{R1}).  It follows from various results in \cite{TW2}
that $A$ is $\mathcal{Z}$-stable, whence $\mathrm{gr}(A)=1$. 

\hspace{3mm}

\noindent
{\bf gr = 2.}
Let $A$ be a simple, unital AH algebra given as the limit of an inductive 
system
\[
\left(p_i(\mathrm{C}(X_i) \otimes \mathcal{K})p_i,\phi_i \right),
\]
where the $X_i$ are compact connected Hausdorff spaces, $p_i \in \mathrm{C}(X_i) 
\otimes \mathcal{K}$ is a projection, and 
\[
\phi_i:p_i(\mathrm{C}(X_i) \otimes \mathcal{K})p_i \to 
p_{i+1}(\mathrm{C}(X_{i+1}) \otimes \mathcal{K})p_{i+1}
\]
is an unital $*$-homomorphism.  Suppose that
\[
\frac{\mathrm{dim}(X_i)}{\mathrm{dim}(p_i)} \stackrel{i \to \infty}{\longrightarrow}
c \in \mathbb{R}, \ c \neq 0.
\]
Since $\mathrm{dim}(p_i) \to \infty$ as $i \to \infty$ by simplicity, we have that
$2 \geq \mathrm{tdg}(A) \geq \mathrm{gr}(A)$.  If $A$ is not $\mathcal{Z}$-stable,
then $\mathrm{gr}(A) =2$.  Many of the known examples of non-$\mathcal{Z}$-stable
simple, nuclear $C^*$-algebras have this form, including the AH algebras of 
\cite{V1} having perforated ordered $\mathrm{K}_0$-groups, those of \cite{V2} having
finite non-minimal stable rank, the algebra $B$ of \cite{R5} which is not stable but
for which $\mathrm{M}_2(B)$ is stable, and the counterexample to Elliott's classification
conjecture in \cite{To2}.  

Let $A$ be a simple, nuclear $C^*$-algebra containing a finite
and an infinite projection and satisfying the UCT (the existence of such algebras is
established in \cite{R3}).  Kirchberg proves in \cite{K} that
the tensor product of any two simple, unital and infinite-dimensional
$C^*$-algebras is either stably finite or purely infinite.  It follows
that $A \otimes A$ is purely infinite and hence $\mathcal{Z}$-stable, 
so $\mathrm{gr}(A) = 2$.

\hspace{3mm}

\noindent
{$\mathbf{gr > 2.}$}  The algebras in Theorem \ref{range} are the first examples of simple nuclear
algebras with finite growth rank strictly greater than 2. 
The algebra of \cite{V2} having infinite stable rank probably also has infinite growth
rank --- it bears a more than passing resemblance to the algebra of infinite growth rank in Theorem \ref{range}.

\section{Tensor factorisation}\label{tens}

A simple $C^*$-algebra is said to be {\it tensorially prime} if it cannot be written
as a tensor product $C \otimes D$, where both $C$ and $D$ are simple and non-type-I.
It has been surprising to find that the majority of our stock-in-trade simple, separable, and
nuclear $C^*$-algebras are not tensorially prime --- every class of simple, separable, and nuclear
$C^*$-algebras for which the Elliott conjecture is currently confirmed consists of 
$\mathcal{Z}$-stable members (\cite{JS1}, \cite{TW2}).  Kirchberg (\cite{K}) has
shown that every simple exact $C^*$-algebra which is not tensorially prime is either
stably finite or purely infinite.  R{\o}rdam has produced an example of a simple
nuclear $C^*$-algebra containing both a finite and an infinite projection which, in
light of Kirchberg's result, is tensorially prime (\cite{R3}).  The question of whether
every simple, nuclear, and non-$\mathcal{Z}$-stable $C^*$-algebra is tensorially prime
has remained open.  Theorem \ref{range} settles this question, negatively.

\begin{cors}\label{nottensprime}
There exists a simple, nuclear, and non-$\mathcal{Z}$-stable $C^*$-algebra 
which is not tensorially prime.  
\end{cors}

\begin{proof}
Let $A$ be the algebra of growth rank three in Theorem \ref{range}.  $A \otimes A$
satisfies the hypotheses of the corollary, yet is evidently not tensorially prime.
\end{proof}

It is interesting to ask whether simple, unital, and nuclear $C^*$-algebras which fail 
to be tensorially prime must in fact have an infinite factorisation, i.e., can
be written as $\otimes_{i=1}^{\infty} C_i$, where each $C_i$ is simple, unital,
nuclear, and non-type-I.  This is trivially true for $\mathcal{Z}$-stable algebras,
since $\mathcal{Z} \cong \mathcal{Z}^{\otimes \infty}$ (cf. \cite{JS1}).  R{\o}rdam
has asked whether every simple, unital, nuclear, and non-type-I $C^*$-algebra admits
an unital embedding of $\mathcal{Z}$.  If this turns out to be true for separable
algebras, then Theorem 7.2.2 of \cite{R1} implies that infinite tensor products
of such algebras are always $\mathcal{Z}$-stable.  This, in turn, will imply that simple, unital,
separable, and nuclear $C^*$-algebras which do not absorb $\mathcal{Z}$ tensorially cannot
have an infinite tensor factorisation. 

\section{Universal properties and infinite tensor products}\label{uniprops}

Little is known about the extent to which $\mathcal{Z}$ is unique, save that it
is determined by its $\mathrm{K}$-theory inside a small class of $\mathcal{Z}$-stable
inductive limit algebras (\cite{M}, \cite{TW2}).  The Elliott conjecture, which may well
hold for the class of simple, separable, nuclear, infinite-dimensional and $\mathcal{Z}$-stable 
$C^*$-algebras, predicts that $\mathcal{Z}$ will be the unique such algebra which is furthermore
unital, projectionless, unique trace, and $\mathrm{KK}$-equivalent to $\mathbb{C}$.

R{\o}rdam has suggested the following universal property, which could conceivably be verified
for $\mathcal{Z}$ within the class of separable, unital, and
nuclear $C^*$-algebras having no finite-dimensional representations.

\begin{uprops}\label{ror} Let $\mathcal{C}$ be a class of separable,
    unital, and nuclear $C^*$-algebras.  
If $A$ in $\mathcal{C}$ is such that
\begin{enumerate}
\item[(i)] every unital endomorphism of $A$ is approximately inner, and
\item[(ii)] every $B$ in $\mathcal{C}$ admits an unital embedding
$\iota: A \to B$,
\end{enumerate}
then $A$ is unique up to $*$-isomorphism.
\end{uprops}

\begin{proof}
Elliott's Intertwining Argument (cf. \cite{EG}).
\end{proof}

We propose a second property.

\begin{uprops}\label{toms}  Let $\mathcal{C}$ be a class of 
unital and nuclear $C^*$-algebras.  If $A$ in $\mathcal{C}$ is such that
\begin{enumerate}
\item[(i)] $A^{\otimes \infty} \cong A$, and
\item[(ii)] $B^{\otimes \infty} \otimes A \cong B^{\otimes \infty}$ for every $B$ in
$\mathcal{C}$,
\end{enumerate}
then $A$ is unique up to $*$-isomorphism.
\end{uprops}

\begin{proof}
Suppose that $A,B$ in $\mathcal{C}$ satisfy (i) and (ii) above.  Then,
\[
A \stackrel{\mathrm{(i)}}{\cong} A^{\otimes \infty} \stackrel{\mathrm{(ii)}}{\cong}
A^{\otimes \infty} \otimes B \stackrel{\mathrm{(i)}}{\cong} A \otimes
B^{\otimes \infty} \stackrel{\mathrm{(ii)}}{\cong} B^{\otimes \infty}
\stackrel{\mathrm{(i)}}{\cong} B.
\]
\end{proof}

Universal Properties \ref{ror} and \ref{toms} have the same basic
structure.  In each case, condition (i) is intrinsic and known to 
hold for the Jiang-Su algebra $\mathcal{Z}$, while condition (ii) is
extrinsic and potentially verifiable for $\mathcal{Z}$.
Conditions \ref{ror} (i) and \ref{toms} (i) are skew, but not completely
so:  the
separable, unital $C^{*}$-algebras satisfying both conditions are
precisely the strongly self-absorbing $C^{*}$-algebras studied in
\cite{TW1}.  Any separable, unital, and nuclear $C^{*}$-algebra $A$ which
admits an unital embedding of $\mathcal{Z}$ then satisfies $A^{\otimes 
\infty} \otimes \mathcal{Z} \cong A^{\otimes \infty}$ (cf.  Theorem
7.2.2 of Rordam);  if Universal Property \ref{ror} is satisfied by $\mathcal{Z}$
inside a class $\mathcal{C}$ of separable, unital, and nuclear
$C^{*}$-algebras, then the same is true of
Universal Property \ref{toms}.  The attraction of condition \ref{toms} (ii), as we shall see, is that it can be
verified (with $A = \mathcal{Z}$) for a large class of projectionless $C^{*}$-algebras;
there is, to date, no similar confirmation of condition \ref{ror} (ii).

An interesting point:  if one takes $\mathcal{C}$ to be
the class of Kirchberg algebras, then Universal Properties \ref{ror} and
\ref{toms} both identify $\mathcal{O}_{\infty}$;  if one takes
$\mathcal{C}$ to be the class of simple, nuclear, separable, and
unital $C^{*}$-algebras satisfying the Universal Coefficients Theorem 
and containing an infinite projection --- a class which properly contains
the Kirchberg algebras --- then Universal Property \ref{toms}
still identifies $\mathcal{O}_{\infty}$, while Universal Property
\ref{ror} does not (indeed, it is unclear whether \ref{ror} identifies
anything at all in this case).

To prove that $\mathcal{Z}$ satisfies Universal
Property \ref{toms} among unital and nuclear $C^{*}$-algebras,
one must determine whether infinite tensor products of such algebras
are $\mathcal{Z}$-stable.  Formally, the question is reasonable.  If $\mathrm{gr}(A^{\otimes \mathrm{gr}(A)}) 
= 1$ whenever $\mathrm{gr}(A) < \infty$, then why not $\mathrm{gr}(A^{\otimes \infty}) = 1$ 
for any $A$?  It follows immediately from Definition \ref{growthrank} that
one has either $\mathrm{gr}(A^{\otimes \infty}) = 0$ or $\mathrm{gr}(A^{\otimes \infty})=\infty$.  

Recall that for natural numbers $p,q,n$ such that $p$ and $q$ divide $n$, the \emph{dimension drop interval} 
$\mathrm{I}[p,n,q]$ is the algebra of functions
\[
\{ f \in \mathrm{C}([0,1],\mathrm{M}_n)|f(0) = a \otimes 1_{n/q}, a \in \mathrm{M}_p, f(1) = b \otimes 1_{n/p},
b \in \mathrm{M}_q \}.
\]
If $p$ and $q$ are relatively prime and $n = pq$, then we say that $\mathrm{I}[p,pq,q]$ is
a \emph{prime dimension drop interval}.  $\mathcal{Z}$ is the unique simple and unital 
inductive limit of prime dimension drop intervals having 
\[
(\mathrm{K}_0,\mathrm{K}_0^+,[1]) \cong (\mathbb{Z},\mathbb{Z}^+,1); \ \ \ \mathrm{K}_1 = 0; \ \ \
\mathrm{T} = \{*\}
\]
(\cite{JS1}).

The next two propositions are germane to the results in this section.  They follow more 
or less directly from Proposition 2.2 and Theorem 2.3 of \cite{TW2}, respectively.

\begin{props}\label{Zembed-grinf=0}
Let $A$ be a separable, nuclear, and unital $C^*$-algebra.  Then,
$\mathrm{gr}(A^{\otimes \infty})=1$ if and only if there exists, for any relatively 
prime natural numbers $p$ and $q$, an unital 
$*$-homomorphism $\iota:\mathrm{I}[p,pq,q] \to 
A^{\otimes \infty}$. 
\end{props}

\begin{proof}
The ``only if'' part of the proposition is straightforward --- every prime dimension
drop interval embeds into $\mathcal{Z}$, which in turn embeds into $A^{\otimes \infty} \otimes
\mathcal{Z} \cong A^{\otimes \infty}$.

Proposition 2.2 of \cite{TW2} states:  if $B$ is a separable and unital $C^*$-algebra and there
exists, for each pair of relatively prime natural numbers $p$ and $q$, an unital $*$-homomorphism
\[
\phi: I[p,pq,q] \to \frac{\prod_{i=1}^{\infty} B}{\bigoplus_{i=1}^{\infty} B} \cap B',
\]
where $B'$ is the commutant of embedding of $B \to \prod_{i=1}^{\infty} B/\oplus_{i=1}^{\infty} B$
coming from constant sequences, then $B \otimes \mathcal{Z} \cong B$.  
Equivalently, if one has finite sets $F_1 \subseteq F_2
\subseteq \cdots \subseteq B$
such that $\cup_i F_i$ is dense in $B$ and, for any
relatively prime natural numbers $p$ and $q$ and $i \in \mathbb{N}$ 
an unital $*$-homomorphisms $\phi_{p,q} : I[p,pq,q] \to B$ such that 
$\mathrm{Im}(\phi_{p,q})$ commutes with $F_i$ up to $1/2^i$, then $B \otimes \mathcal{Z}
\cong B$.  

Put $B = A^{\otimes \infty}$, and choose finite sets $F_i \subseteq A^{\otimes \infty}$
with dense union.  We may write
\[
A^{\otimes \infty} \cong A^{\otimes \infty} \otimes A^{\otimes \infty} \otimes \cdots,
\]
and assume that that $F_i$ is contained in the first $i$ tensor factors of $A^{\otimes \infty}$
above.
By assumption there exists, for any relatively prime natural numbers $p$ and $q$, an unital
$*$-homomorphism $\phi: I[p,pq,q] \to A^{\otimes \infty}$.  By composing $\phi$ with the embedding
of $A^{\otimes \infty}$ as the $(i+1)^{\mathrm{th}}$ tensor factor of $A^{\otimes \infty} 
\otimes A^{\otimes \infty} \otimes \cdots$, we obtain an unital $*$-homomorphism from
$I[p,pq,q]$ to $A^{\otimes \infty}$ whose image commutes with $F_i$, as required.  Thus,
$A^{\otimes \infty} \otimes \mathcal{Z} \cong A^{\otimes \infty}$ and $\mathrm{gr}(A^{\otimes 
\infty}) = 1$.
\end{proof}

\begin{props}\label{fd-grproduct-zero} 
Let $A$ be a separable, nuclear, and unital $C^*$-algebra.  Suppose
that $A$ admits an unital $*$-homomorphism $\iota:\mathrm{M}_{2} \oplus
\mathrm{M}_{3} \to A$.  Then, $\mathrm{gr}(A^{\otimes \infty}) = 1$.
\end{props}

\begin{proof}
Choose finite sets $F_1 \subseteq F_2 \subseteq \cdots \subseteq A^{\otimes \infty}$
with dense union, and with the property that $F_i$ is contained in the first $i$ tensor
factors of $A^{\otimes \infty}$.  One can then use $\iota$ to obtain an unital 
$*$-homomorphism from $\mathrm{M}_2 \oplus \mathrm{M}_3$ to the $(i+1)^{\mathrm{th}}$
tensor factor of $\mathrm{A^{\otimes \infty}}$.  In particular, the image of this
homomorphism commutes with $F_i$.  It follows that $A^{\otimes \infty}$
is approximately divisible, and hence $\mathcal{Z}$-stable by Theorem 2.3 of \cite{TW1}.
\end{proof}

\begin{cors}\label{perror}
Let $A$ be a separable, nuclear and unital $C^*$-algebra of real rank zero 
having no one-dimensional representation.  Then, $\mathrm{gr}(A^{\otimes \infty})=1$.
\end{cors}

\begin{proof}
In Proposition 5.7 of \cite{RP}, Perera and R{\o}rdam prove that an algebra
$A$ as in the hypotheses of the corollary admits an unital embedding of a
finite-dimensional algebra $F$ having no direct summand of dimension one.
Apply Proposition \ref{fd-grproduct-zero}.
\end{proof}

It it not known at present whether every simple and unital AH algebra admits an unital
embedding of $\mathcal{Z}$.  We will prove that infinite tensor products of such algebras
are nevertheless $\mathcal{Z}$-stable whenever they lack one-dimensional
representations.

\begin{lms}\label{homembed}
Given any natural number $N$, there exists $\epsilon>0$ with the following 
property:  if
\[
A:=p(\mathrm{C}(X) \otimes \mathcal{K})p,
\]
$X$ a connected finite CW-complex, is such that 
\[
\frac{\mathrm{dim}(X)}{\mathrm{rank}(p)} < \epsilon,
\]
then there is an unital $*$-homomorphism
\[
\iota: \mathrm{M}_{N} \oplus \mathrm{M}_{N+1} \to A.
\]
\end{lms}

\begin{proof}
Since $X$ is a finite CW-complex, 
the $\mathrm{K}_{0}$-group of $A$ is finitely generated.  Write
\[
\mathrm{K}_{0}A \cong G_{1} \oplus G_{2} \oplus \cdots \oplus G_{k}, 
\]
where each $G_{i}$ is cyclic, and $G_{1} = \langle [\theta_{1}] 
\rangle$ is the free cyclic group
generated by the $\mathrm{K}_{0}$-class $[\theta_{1}]$
of the trivial complex line bundle $\theta_{1}$ over $X$.  

Let $\mathrm{rank}(p)$ be large enough --- equivalently, 
$\mathrm{dim}(X)/\mathrm{rank}(p)$ small enough --- to
ensure the existence of non-negative integers $a$ and $b$ such that
\[
\mathrm{rank}(p)  =  aN + b(N+1), \
a,b \geq \mathrm{dim}(X)/2.
\]
Write
\[
[p] = \oplus_{j=1}^{k} g_{j}, \ g_{j} \in G_{j}, \ 1 \leq j \leq m.
\]
One has, by definition, that $g_{1}= Na[\theta_{1}] + (N+1)b[\theta_{1}]$.
Since $N$ and $N+1$ are relatively prime one also has, for every $i \geq
2$, elements $h_{i}, r_{i}$ of $G_{i}$ such that 
\[
g_{i} = Nh_{i} + (N+1)r_{i}.
\]
Set
\[
h := h_{1} \oplus \cdots \oplus h_{k}; \ \ r := r_{1} \oplus \cdots
\oplus r_{k}.
\]
Then, $g = Nh +(N+1)r$, and $h, r \in \mathrm{K}_{0}(A)^{+}$ --- 
the virtual dimension of these elements exceeds $\mathrm{dim}(X)/2$.

Find pairwise orthogonal projections $P_{1},\ldots,P_{N}$ in
$\mathrm{M}_{\infty}(A)$ such that $[P_{i}]=h$, $1 \leq i \leq N$.
Similarly, find pairwise orthogonal projections $Q_{1},\ldots,Q_{N+1}$
such that $[Q_{j}]=r$, $1 \leq j \leq N+1$.  Since
$(\oplus_{i} P_{i}) \oplus (\oplus_{j} Q_{j})$ is Murray von-Neumann
equivalent to $p$, we may assume that the $P_{i}$s and $Q_{j}$s are
in $A$.  Furthermore, $P_{i}$ and $P_{k}$ are Murray-von Neumann
equivalent for any $i$ and $k$, and similarly for $R_{i}$ and $R_{k}$.
One may then easily find a system of matrix units for $\mathrm{M}_{N}$
and $\mathrm{M}_{N+1}$ using the partial isometries implementing the
equivalences among the $P_{i}$s and $R_{j}$s.  It follows that
there is an unital embedding of $\mathrm{M}_{N} \oplus \mathrm{M}_{N+1}$
into $A$. 
\end{proof}

\begin{props}\label{lochom}
Let $A$ be a separable, unital $C^{*}$-algebra.   
Let  
\[
B = \bigoplus_{i=1}^{n} p_{i}(\mathrm{C}(X_{i}) \otimes
\mathcal{K})p_{i}
\]
satisfy $\mathrm{rank}(p_{i}) \geq 2$.  If there is an unital $*$-homomorphism $
\phi: B \to A$, then $\mathrm{gr}(A^{\otimes \infty}) = 1$.
\end{props}

\begin{proof}
For any natural number $k$ one has an unital $*$-homomorphism
\[
\phi^{\otimes k}: B^{\otimes k} \to A^{\otimes k}.
\]
Let $\iota_{k}: A^{\otimes k} \to A^{\otimes \infty}$ be the map 
obtained by embedding $A^{\otimes k}$ as the first $k$ factors of 
$A^{\otimes \infty}$.  Setting $\gamma_{k}=\iota_{k} \circ
\phi^{\otimes k}$, one has an unital $*$-homomorphism 
\[
\gamma_{k}: B^{\otimes k} \to A^{\otimes \infty}.
\]

Recall that
\begin{equation}\label{homprod}
p(\mathrm{C}(X) \otimes \mathcal{K})p \otimes q(\mathrm{C}(Y) \otimes 
\mathcal{K})q \cong (p \otimes q)(\mathrm{C}(X \times Y) \otimes
\mathcal{K})(p \otimes q)
\end{equation}
for compact Hausdorff spaces $X$ and $Y$ and projections $p
\in \mathrm{C}(X) \otimes \mathcal{K}$, $q \in \mathrm{C}(Y)
\otimes \mathcal{K}$.  Let $Z$ be any connected component of the 
spectrum of $B^{\otimes k}$, and let $p_{Z} \in B^{\otimes k}$ be
the projection which is equal to the unit of $B^{\otimes k}$ at every point in $Z$, and
equal to zero at every other point in the spectrum of $B^{\otimes k}$.
It follows from equation \ref{homprod} that
\[
\frac{\mathrm{dim}(Z)}{\mathrm{rank}(p_{Z})} \leq
\frac{k \left(\mathrm{max}_{1 \leq i \leq n}
\mathrm{dim}(X_{i})\right)}{\left(\mathrm{min}_{1 \leq i \leq n}
\mathrm{rank}(p_{i})\right)^{k}} \leq
\frac{k \left(\mathrm{max}_{1 \leq i \leq n}
\mathrm{dim}(X_{i})\right)}{2^{k}} \stackrel{k \to
\infty}{\longrightarrow} 0.
\]
Thus, for a fixed $N \in \mathbb{N}$, there is some $k \in \mathbb{N}$ 
such that every homogeneous direct summand of $B^{\otimes k}$ with
connected spectrum satisfies 
the hypothesis of Lemma \ref{homembed} for the corresponding value of
$\epsilon$.  It follows that there is an unital $*$-homomorphism 
\[
\psi: \mathrm{M}_{N} \oplus \mathrm{M}_{N+1} \to B^{\otimes k}.
\]
The composition $\gamma_{k} \circ \psi$ yields an unital
$*$-homomorphism from $\mathrm{M}_{N} \oplus \mathrm{M}_{N+1}$ to
$A^{\otimes \infty}$  (we may assume that $N \geq 2$).  It follows
that there is an unital $*$-homomorphism from $\mathrm{M}_2 \oplus \mathrm{M}_3$
to $\mathrm{M}_N \oplus \mathrm{M}_{N+1}$, and hence an unital $*$-homomorphism
from $\mathrm{M}_2 \oplus \mathrm{M}_3$ to $A^{\otimes \infty}$.
Apply Proposition \ref{fd-grproduct-zero} to conclude that
$(A^{\otimes \infty})^{\otimes \infty} \cong A^{\otimes \infty}$ is 
$\mathcal{Z}$-stable.
\end{proof}

An algebra $B$ as in the statement of Proposition \ref{lochom} need not
have any non-trivial projections.  Take, for instance, the algebra
$p(\mathrm{C}(S^4) \otimes \mathcal{K})p$, where $p$ is the higher Bott projection;  $p$
has no non-zero Whitney summands by a Chern class argument.
On the other hand, the proof of Proposition \ref{lochom} shows that if
$A$ satisfies the hypotheses of the same, then $A^{\otimes \infty}$
has many projections.

\begin{thms}\label{AHnooned}
Let $A$ be an unital AH algebra having no one-dimensional representation.
Then, $\mathrm{gr}(A^{\otimes \infty})=1$.
\end{thms}

\begin{proof}
Write
\[
A = \lim_{i \to \infty} \left( A_{i} := \bigoplus_{l=1}^{n_{i}}
p_{i,l}(\mathrm{C}(X_{i,l}) \otimes \mathcal{K})p_{i,l}, \ \phi_{i}
\right),
\]
where $\phi_i:A_i \to A_{i+1}$ is unital.
Define $J_{i}:=\{l \in \mathbb{N}| \mathrm{rank}(p_{i,l})=1\}$, and put 
\[
B_{i} := \bigoplus_{l \in J_{i}} p_{i,l}(\mathrm{C}(X_{i,l}) \otimes
\mathcal{K})p_{i,l}.
\]
Define $\psi_{i}: B_{i} \to B_{i+1}$ by restricting $\phi_i$ to $B_i$, then
cutting down the image by the unit of $B_{i+1}$.  Put $\psi_{i,j}:=
\psi_{j-1} \circ \cdots \circ \psi_{i}$.  Notice that for reasons of rank, each $\psi_i$ is
unital --- the only summands of $A_i$ whose images in $B_{i+1}$ may
be non-zero are the summands which already lie in $B_{i}$.

If, for some $i \in \mathbb{N}$, one has $\psi_{i,j} \neq 0$ for
every $j>i$, then one may find, for every $j>i$, a rank one 
projection $q_{j} \in \{p_{j,l}\}_{l=1}^{n_{j}}$ such that 
the cut-down of the image of $\psi_{j}|_{q_{j}B_{j}q_{j}}$
by $q_{j+1}$ gives an unital $*$-homomorphism from 
$q_{j}B_{j}q_{j}$ to $q_{j+1}B_{j+1}q_{j+1}$.
Let $Y_{j}$ be the spectrum of $q_{j}B_{j}q_{j}$.  There is a
continuous map $\theta_{j}:Y_{j+1} \to Y_{j}$ such that 
$\psi_{j}(f)(y) = f(\theta_{j}(y))$ for every $y \in Y_{j+1}$
and $f \in q_j B_j q_j$.
Choose a sequence of points $y_{j} \in Y_{j}$, $j>i$, with the
property that $\theta_{j}(y_{j+1}) = y_{j}$.  Let $\gamma_{j}: 
A_{j} \to \mathbb{C}$ be given by $\gamma_{j}(f) = f(y_{j})$.
Then $(\gamma_{j})_{j>i}$ defines an unital inductive limit $*$-homomorphism
$\gamma: A \to \mathbb{C}$; $A$ has a one-dimensional representation.
We conclude that for every $i \in \mathbb{N}$ there exists $j > i$
such that $\psi_{i,j} = 0$.  It follows that $B_{j} = \{0\}$, so
that $A_{j}$ has no one-dimensional representations.   Apply Proposition \ref{lochom}
to conclude that $\mathrm{gr}(A^{\otimes \infty})=1$.
\end{proof}

\noindent
Theorem \ref{AHnooned} is interesting in light of the fact
that there are unital AH algebras which are not weakly divisible
(every algebra constructed in \cite{V2} has this deficiency, for instance), 
so Proposition \ref{fd-grproduct-zero} cannot be applied to them.

In the case of simple, unital AH algebras, infinite tensor products
are not only $\mathcal{Z}$-stable, but classifiable as well.  The next
proposition has been observed independently by Bruce Blackadar and the
author.

\begin{props}\label{AHclass}
Let $A$ be a simple, unital AH algebra.  Then, $A^{\otimes \infty}$
has very slow dimension growth in the sense of \cite{G}.  
\end{props}

\begin{proof}  Write $A = \lim_{i \to \infty} (A_i,\phi_i)$ where, as usual,
\[
A_i = \bigoplus_{i=1}^{m_i} p_{i,l}(\mathrm{C}(X_{i,l}) \otimes \mathcal{K})p_{i,l}.
\] 
Define
\[
n_i := \mathrm{min}_{1 \leq l \leq m_i}\{\mathrm{rank}(p_{i,l})\}; \ \ 
k_i := \mathrm{max}_{1 \leq l \leq m_i}\{\mathrm{dim}(X_{i,l})\}.
\]
Let $\epsilon_1, \epsilon_2,\ldots$ be a sequence of positive tolerances converging to zero.
Set $r_1 = 1$, and choose for each $i \in \mathbb{N}$ a natural number $r_i \in \mathbb{N}$ 
satisfying 
\[
\frac{n_i^{r_i}}{(k_i r_i)^3} < \epsilon_i.
\]
$A^{\otimes \infty}$ can be decomposed as follows: 
\[
A^{\otimes \infty} = \lim_{i \to \infty}(A_i^{\otimes r_i}, \phi_i^{\otimes r_i} \otimes 1_{A_{i+1}}^{\otimes r_{i+1} - r_i}).
\]
The maximum dimension of a connected component of the spectrum of $A_i^{\otimes r_i}$ is $r_i k_i$,
while the minimum rank of the unit of a homogeneous direct summand of $A_i^{\otimes r_i}$ corresponding
to such a connected component is $n_i^{r_i}$.  It follows that the decomposition for $A^{\otimes \infty}$
above has very slow dimension growth, whence the limit is approximately divisible (\cite{EGL2}), 
$\mathcal{Z}$-stable (Theorem 2.3 of \cite{TW2}), and classifiable via the Elliott invariant
(\cite{EGL}). 
\end{proof}

Certain infinite tensor products of
$C^{*}$-algebras of real rank zero are also classifiable.  The 
proposition below is direct consequence of recent work by 
Brown (\cite{B}).

\begin{props}\label{prodclass}
Let $A$ be a simple, unital, inductive limit of type-I $C^*$-algebras 
with a unique tracial state.  If there is an unital $*$-homomorphism
$\phi: \mathrm{M}_2 \oplus \mathrm{M}_3 \to A$, then $A^{\otimes \infty}$ is tracially AF.
\end{props}

\begin{proof}
$A^{\otimes \infty}$ is $\mathcal{Z}$-stable by an application of Proposition
\ref{fd-grproduct-zero}, whence it has weakly unperforated $\mathrm{K}_0$-group
(Theorem 1 of \cite{GJS}) and stable rank one (Theorem 6.7 of \cite{R4}).
Since the proof of Proposition \ref{fd-grproduct-zero} actually shows that
$A^{\otimes \infty}$ is approximately divisible, we conclude that it has
real rank zero by the main theorem of \cite{BKR}.  We have thus collected the
hypotheses of Corollary 7.9 of \cite{B}, whence $A^{\otimes \infty}$ is 
tracially AF.
\end{proof}

\noindent
Notice that algebras satisfying the hypotheses of Proposition
\ref{prodclass} need not be tracially AF, even if one excludes
the trite example of a finite-dimensional algebra with
no one-dimensional representation.  Examples include
$\mathrm{M}_{2}(A)$ for
any of the algebras produced in Theorem \ref{range} or any
algebra constructed in \cite{V2}, hence algebras of arbitrary 
growth rank or stable rank.

The infinite tensor products considered so far have all contained
non-trivial projections.  We turn now to certain potentially
projectionless infinite tensor products.   

\begin{dfs}\label{dimdrop}
Let there be given a homogeneous $C^{*}$-algebra
$\mathrm{M}_{k}(\mathrm{C}(X))$ ($X$ is not necessarily
connected) and closed pairwise disjoint sets $X_{1},\ldots,X_{n}
\subseteq X$.  Let $F$ be a finite-dimensional
$C^{*}$-algebra, and let $\iota_{i}:F \to \mathrm{M}_{k}$,
$1 \leq i \leq n$, be unital $*$-homomorphisms.  Define
\[
\phi_{i}:F \to \mathrm{M}_{k}(\mathrm{C}(X_{i})) \cong \mathrm{C}(X_{i})
\otimes \mathrm{M}_{k}
\]
by $\phi_{i}:=\mathbf{1} \otimes \iota_{i}$, and put
\[
A := \{f \in \mathrm{M}_{k}(\mathrm{C}(X))|f|_{X_{i}} \in
\mathrm{Im}(\phi_{i})\}.
\]
We call $A$ a \emph{generalised dimension drop algebra}.
\end{dfs}

Separable and unital direct limits of direct sums of generalised
dimension drop algebras are, in general, beyond the
scope of current methods for classifying approximately subhomogeneous
(ASH) algebras via
$\mathrm{K}$-theory.  The only such algebras which are known
to admit an unital embedding of $\mathcal{Z}$ are those for
which classification theorems exist, and these form a quite limited
class.  But for infinite tensor products, we can prove the following:

\begin{thms}\label{gendropinf}
Let $A$ be a separable, unital, and nuclear $C^{*}$-algebra.  Suppose 
that for every $n \in \mathbb{N}$ there is a finite direct sum of 
generalised dimension drop algebras $B_{n}$ having no representation
of dimension less than $n$, and an unital $*$-homomorphism
$\gamma_{n}:B_{n} \to A$.  Then, $A^{\otimes \infty} \otimes \mathcal{Z} \cong
A^{\otimes \infty}$.
\end{thms}

Theorem \ref{gendropinf} follows directly from Proposition \ref{Zembed-grinf=0}
and Lemma \ref{dimdropembed} below.

\begin{lms}\label{dimdropembed}
Let $I[p,pq,q]$ be a fixed prime dimension drop interval, and let $B$ be 
a generalised dimension drop algebra.  There exists $N \in \mathbb{N}$
such that if every non-zero finite-dimensional
representation of $B$ has dimension at least $N$, then there is a
unital $*$-homomorphism $\gamma:I[p,pq,q] \to B$.
\end{lms}

\begin{proof}
Let $N \geq pq-p-q$.  It is well known that for any natural number
$M \geq N$ there are non-negative integers $a_{M}$ and $b_{M}$ such 
that $a_{M}p + b_{M}q = M$.  Since $B$ has no representations of 
dimension less than $N$, we may assume that every simple direct summand of
the finite-dimensional algebra $F$ associated to $B$ has matrix size
at least $N$.  There is an unital $*$-homomorphism
$\psi: I[p,pq,q] \to F$ defined as follows:  given 
a simple direct summand $\mathrm{M}_{k_{j}}$ of F, $1 \leq j \leq m$,
define a map
\[
\psi_{j}: I[p,pq,q] \to \mathrm{M}_{k_{j}}
\]
by 
\[
\psi_{j}(f) = \bigoplus_{l=1}^{a_{k_{j}}} f(0) \oplus
\bigoplus_{r=1}^{b_{k_{j}}} f(1);
\]
put
\[
\psi := \bigoplus_{j=1}^{m} \psi_{j}.
\]

Adopt the notation of Definition \ref{dimdrop} for $B$.
Find pairwise disjoint open sets $O_{i} \supseteq X_{i}$, $1 \leq i
\leq n$, and put $C = (\cup_{i} O_{i})^{c}$.  Since $X$ is normal, 
there is a continuous function $f:X \to [0,1]$ such that $f = 0$
on $C$ and $f = 1$ on $\cup_{i} X_{i}$.  Define a map 
\[
\gamma_{1}: I[p,pq,q] \to \mathrm{M}_{k}(\mathrm{C}(C \cup O_{1}))
\]
by 
\[
\gamma_{1}:= (\mathbf{1}_{\mathrm{M}_{k}(\mathrm{C}(C \cup O_{1}))} \otimes 
\iota_{1}) \circ \psi.
\]
For each $2 \leq i \leq n$ define similar maps
\[
\gamma_{i}: I[p,pq,q] \to \mathrm{M}_{k}(\mathrm{C}(X_{i}))
\]
by
\[
\gamma_{i}:= (\mathbf{1}_{\mathrm{M}_{k}(\mathrm{C}(X_{i}))} \otimes 
\iota_{i}) \circ \psi.
\]
Lemma 2.3 of  \cite{JS1} shows that the space of unital
$*$-homomorphisms from $I[p,pq,q]$ to $\mathrm{M}_{k}$ is contractible.
Choose, then, for each $2 \leq i \leq n$, a homotopy
\[
\omega_{k}: [0,1] \times I[p,pq,q] \to \mathrm{M}_{k}
\]
such that $\omega_{i}(0,g) = \gamma_{1}(g)$ and $\omega_{i}(1,g) =
\gamma_{i}(g)$, $\forall g \in I[p,pq,q]$.  Define an unital
$*$-homomorphism $\gamma: I[p,pq,q] \to B$ by
\[
\gamma(g)(x) = \left\{ 
\begin{array}{rl} 
\gamma_{1}(g), & \ x \in C \cup O_{1} \\
\gamma_{i}(g), & \ x \in X_{i}, \ 2 \leq i \leq n \\
\omega_{i}(t,g), & \ x \in O_{i} \backslash X_{i} \ \mathrm{and} \
f(x) = t
\end{array}
\right.
\]
\end{proof}

The hypotheses of Theorem \ref{gendropinf} are less general than one 
would like --- replacing generalised dimension drop algebras with 
recursive subhomogeneous algebras would be a marked improvement.  
On the other hand, they do not even require that the algebra $A$ be approximated
locally on finite sets by generalised dimension drop algebras, and are satisfied
by a wide range of $C^*$-algebras:
\begin{enumerate}
\item[(i)]
simple and unital limits of inductive sequences $(A_i, \phi_i)$, $i \in \mathbb{N}$,
where each $A_i$ is a finite direct sum of generalised dimension drop algebras 
--- these encompass all approximately subhomogeneous (ASH) algebras
for which the Elliott conjecture is confirmed;
\item[(ii)] for every weakly unperforated instance of the Elliott invariant $I$, an unital,
separable, and nuclear $C^*$-algebra $A_I$ having this invariant (see the proof of the main
theorem of section 7 in \cite{EV});
\item[(iii)] simple, unital, separable, and nuclear $C^*$-algebras having the same Elliott
invariant as $\mathcal{Z}$ for which there are no classification results (the main theorem
in section 7 of \cite{EV} provides a construction of a simple, unital, separable, and
nuclear $C^*$-algebra with the same Elliott invariant as $\mathcal{Z}$;  there are no ASH
classification results which cover this algebra, yet it satisfies the hypotheses of 
Theorem \ref{gendropinf}).
\end{enumerate}

\section{The growth rank and other ranks}\label{otherranks}

In this last section we explore the connections between the growth rank and other
ranks for nuclear $C^*$-algebras:  the stable rank ($\mathrm{sr}(\bullet)$), 
the real rank ($\mathrm{rr}(\bullet)$), the tracial topological rank 
($\mathrm{tr}(\bullet)$), and the decomposition rank ($\mathrm{dr}(\bullet)$). 

Growth rank one is the condition that $A$ absorbs $\mathcal{Z}$ tensorially.
If, in addition, $A$ simple and unital, then it is either stably finite 
or purely infinite by Theorem 3 of \cite{GJS}.
If $A$ is purely infinite, then it has infinite stable rank and real rank zero
(see \cite{R1}, for instance).  If $A$ is finite, then it has stable rank one 
by Theorem 6.7 of \cite{R4}.  The bound $\mathrm{rr}(A) \leq 2\mathrm{sr}(A) -1$
holds in general, so one also has $\mathrm{rr}(A) \leq 1$ (\cite{BP}).

As mentioned at the end of section \ref{rangesec}, the AH algebras of \cite{V1} with perforated ordered
$\mathrm{K}_0$-groups of bounded perforation and those of \cite{V2} having finite stable rank
all have growth rank two.  For each natural number
$n$ there is an algebra from either \cite{V1} or \cite{V2} with stable rank $n$.
Thus, there is no restriction on the stable rank of a
$C^*$-algebra with growth rank two other than the fact that it is perhaps not infinite. 
The algebra in \cite{V2} of stable rank $n \geq 2$ has real rank equal to $n$ or $n-1$,
so algebras of growth rank two may have more or less arbitrary finite non-zero real
rank.

We can compute the stable rank, real rank, decomposition rank, and tracial 
topological rank of the algebras constructed in Theorem \ref{range}.

\begin{props}\label{infinite-ranks}
Let $A_n$ be the algebra of Theorem \ref{range} such that
$\mathrm{gr}(A) = n \in \mathbb{N} \cup \{\infty\}$.  Then, 
\begin{enumerate}
\item $\mathrm{sr}(A_1) = 1$, $\mathrm{sr}(A_2) < \infty$, and 
$\mathrm{sr}(A_n) = \infty$ for all $n \geq 3$;
\item $\mathrm{rr}(A_1) = 1$, $\mathrm{rr}(A_2) < \infty$, and 
$\mathrm{rr}(A_n) = \infty$ for all $n \geq 3$;
\item $\mathrm{tr}(A_1) = 0$ and $\mathrm{tr}(A_n) = \infty$ for all
$n \geq 2$;
\item $\mathrm{dr}(A_1) = 0$ and $\mathrm{dr}(A_n) = \infty$ for all
$n \geq 2$.
\end{enumerate}
\end{props}

\begin{proof}
Since $A_1$ was taken
to be a UHF algebra, it has stable rank one.  
Inspection of the construction of $A_2$
shows that the ratio 
\[
\frac{\mathrm{dim}(X_i)}{\mathrm{dim}(p_i)}
\]
is bounded above by a constant $K \in \mathbb{R}^+$.  In \cite{N}, Nistor proves that
\[
\mathrm{sr}( p(\mathrm{C}(X) \otimes \mathcal{K})p) = \left\lceil \frac {\lfloor \mathrm{dim}(X)/2 \rfloor}
{\mathrm{rank}(p)} \right\rceil + 1,
\]
where $X$ is a compact Hausdorff space.  It follows that each building
block in the inductive sequence for $A_2$ has stable rank less than $K$,
whence $\mathrm{sr}(A_2) < K < \infty$ by Theorem 5.1 of \cite{Ri1}.
For $n \geq 2$ the algebra $A_n$ is similar to the algebra of infinite stable
rank constructed in Theorem 12 of \cite{V2}.  In fact, the proof of the latter
can be applied directly to show that $\mathrm{sr}(A_n) = \infty$ --- one only 
needs to know that $e(\xi_{I_1^j})^j \neq 0$, which follows from
Proposition 3.2 of \cite{R3} and the fact that $|I_1^j| \geq j$.

The case of the real rank is similar to that of the stable rank.  $A_1$ is UHF,
and so $\mathrm{rr}(A_0) = 0$.  Since the bound $\mathrm{rr}(\bullet) \leq 
2\mathrm{sr}(\bullet) -1$ holds in general, we have $\mathrm{rr}(A_2) < \infty$.
Finally, in a manner analogous to the stable rank case, the proof of Theorem 13 
of \cite{V2} can be applied directly to $A_n$ whenever $n \geq 2$ to show that 
$\mathrm{rr}(A_n) = \infty$.

All UHF algebras have $\mathrm{tr} = 0$, whence $\mathrm{tr}(A_1) = 0$ (\cite{Li3}).
Theorem 6.9 of \cite{Li3} asserts that an unital simple $C^*$-algebra $A$ with
$\mathrm{tr}(A) < \infty$ must have stable rank one.  Since $\mathrm{sr}(A_n) > 1$
for all $n \geq 2$, we conclude that $\mathrm{tr}(A_n) = \infty$ for all such $n$.  
The proof of Theorem 8 of \cite{V2} can be applied directly to $A_2$ to show that
$\mathrm{sr}(A_2) \geq 2$, whence $\mathrm{tr}(A_2) = \infty$. 

All UHF algebras have $\mathrm{dr} = 0$ by Corollary 6.3 of \cite{Wi4}, whence
$\mathrm{dr}(A_1) = 0$.  The same Corollary implies that $\mathrm{dr}(A_n) = \infty$
for all $n \geq 1$ since, by construction, these $A_n$s have unique trace and contain
projections of arbitrarily small trace. 
\end{proof}

Thus, the growth rank is able to distinguish 
between simple $C^*$-algebras which are undifferentiated by other ranks. 
We offer a brief discussion of these other ranks as they relate to Elliott's 
classification program for separable nuclear $C^*$-algebras, and argue that the growth rank meshes 
most naturally with this program.  It must be stressed, however, that these other ranks have
been indispensable to the confirmation of Elliott's conjecture over the years. 

Simple, nuclear, unital and separable $C^*$-algebras of real rank zero have 
so far confirmed Elliott's conjecture, but as the conjecture has also been confirmed for large
classes of simple, nuclear $C^*$-algebras of real rank one, one must conclude that 
real rank zero is at best too strong to characterise the simple, separable and nuclear
$C^*$-algebras satisfying the Elliott conjecture.  There is a counterexample to Elliott's conjecture 
having $\mathrm{rr} = \mathrm{sr} = 1$ in \cite{To2}, so the conditions $\mathrm{rr} \leq 1$ 
and $\mathrm{sr}=1$ also fail to characterise classifiability.

The condition $\mathrm{tr} = 0$ has been shown to be sufficient for the classification
of large swaths of simple, separable, nuclear $C^*$-algebras of real rank zero, but 
the fact that algebras with $\mathrm{tr} < \infty$ must have small projections means 
that this condition will not characterise classifiability --- the Elliott conjecture 
has been confirmed for large classes of projectionless $C^*$-algebras.

The condition $\mathrm{dr} < \infty$ is so far promising as far as characterising
classifiability is concerned.  It may be true that 
\[
\mathrm{dr} < \infty \Leftrightarrow \mathrm{gr} = 1
\]
for simple, separable, nuclear $C^*$-algebras.
On the other hand it remains unclear whether the decomposition rank takes more 
than finitely many values for such algebras, and it is unlikely to distinguish
non-$\mathcal{Z}$-stable algebras.

The growth rank has strong evidence to recommend it as the correct notion of rank
vis a vis classification:  all simple, separable, nuclear and non-elementary
$C^*$-algebras for which the Elliott conjecture is confirmed have $\mathrm{gr}=1$;
all known counterexamples to the conjecture have $\mathrm{gr} > 1$;  the growth
rank achieves every possible value in its range on simple, nuclear and 
separable $C^*$-algebras and is well behaved with respect to common operations.
Furthermore, the classification results available for $\mathcal{Z}$-stable 
$C^*$-algebras are very powerful.  Let $\mathcal{E}$ denote the class of simple,
separable, nuclear, and unital $C^*$-algebras in the bootstrap class $\mathcal{N}$ 
which are $\mathcal{Z}$-stable.  (In light of known examples, $\mathcal{E}$
is the largest class of simple unital algebras for which one can expect the Elliott
conjecture to hold.)  Then:
\begin{enumerate}
\item[(i)] the subclass of $\mathcal{E}_{inf}$ of $\mathcal{E}$ consisting of algebras containing an infinite
projection satisfies the Elliott conjecture (\cite{P1}, \cite{K});
\item[(ii)] the subclass of $\mathcal{E} \backslash 
\mathcal{E}_{inf}$ consisting of algebras with real rank zero and locally
finite decomposition rank satisfies the Elliott conjecture (\cite{Wi6}).
\end{enumerate}

\end{document}